\documentclass[a4paper, final, 12pt]{article}
\usepackage{amsmath, amssymb, latexsym, amscd, amsthm,amsfonts,amstext}
\usepackage[mathscr]{eucal}
\usepackage{graphicx}
\usepackage{subfig}
\usepackage{float}
\usepackage{color}
\usepackage{hyperref}
\usepackage[utf8]{inputenc}
\usepackage[english]{babel}

\numberwithin{equation}{section}
\input{tcilatex}
\begin{document}

\title{A Coefficient Inverse Problem for the Mean Field Games System }
\author{ Michael V. Klibanov \thanks{
Department of Mathematics and Statistics, University of North Carolina at
Charlotte, Charlotte, NC, 28223, USA, mklibanv@charlotte.edu}}
\date{}
\maketitle

\begin{abstract}
A Coefficient Inverse Problem (CIP) of the determination of a coefficient of
the Mean Field Games System (MFGS) of the second order is considered. The
input data are generated by a single measurement event. Lateral Cauchy data,
i.e. Dirichlet and Neumann boundary data are given for solutions of both
equations forming the MFGS. In addition, solutions of both these equations
are given at a fixed moment of time. H\"{o}lder stability estimates are
obtained for both complete and incomplete lateral Cauchy data. These
estimates imply uniqueness of our CIP. Furthermore, our H\"{o}lder stability
estimates can be arranged to become almost Lipschitz stability estimates.
The apparatus of Carleman estimates is the main mathematical tool of this
paper.
\end{abstract}



\section{Introduction}

\label{sec:1}

The mean field games (MFG) theory is a relatively new field, which studies
the collective behavior of large populations of rational decision-makers.
This theory was first introduced in 2006-2007 in seminal publications of
Lasry and Lions \cite{LL3} as well as of Huang, Caines and Malham\'{e} \cite%
{Huang1}. Social sciences enjoy a rapidly increasing role in the modern
society. Therefore, mathematical modeling of social phenomena can
potentially provide a quite important societal impact. In this regard, the
MFG\ theory is the single mathematical model of social processes, which is
based on an universal system of coupled Partial Differential Equations
(PDEs) \cite{Burger}. This is the so-called Mean Field Games system (MFGS).
The number of applications of the MFGS to the societal problems is
flourishing and includes such areas as, e.g. finance, fight with corruption,
cybersecurity, election dynamics, quantum information theory, robotic
control, etc., see, e.g. \cite{A,Burger,Chow,KM,Kol,LL3,Trusov} and
references cited therein.

Thus, due to a broad range of applications of the MFGS, it is important to
address various mathematical questions for this system. One of these
questions is addressed in the current paper. The author addresses in this
paper the question of H\"{o}lder stability estimates and uniqueness theorems
for a Coefficient Inverse Problem (CIP) for the MFGS of the second order.
Stability estimates with respect to possible errors/noise in the input data
are proven here. They can also be considered as accuracy estimates. Such
estimates are important since the input data are always noisy. In addition,
our stability estimates imply uniqueness of our CIP.

Furthermore, our H\"{o}lder stability estimates can be arranged to become
almost Lipschitz stability estimates, see Remark 2.2 in subsection 2.3 as
well as item 3 in Remarks 3.3 in subsection 3.3.

The MFGS of the second order is a system of two coupled nonlinear parabolic
PDEs with two different directions of time \cite{A,LL3}. Let $\Omega \subset 
\mathbb{R}^{n}$ be the domain of interest where agents are located (not
necessary in the physical sense). Let $x\in \Omega $ be the location of an
agent and $t\in \left( 0,T\right) $ be time. The first equation, which is
the Hamilton-Jacobi-Bellman (HJB) equation, describes the behavior of the
price/value function $u\left( x,t\right) $. The second equation, which is
Fokker-Plank equation, describes the density $m\left( x,t\right) $ of agents
in space and time.

The so-called global and local interaction terms depending on $m\left(
x,t\right) $ are parts of the HJB equation. The local interaction term
usually has the form $f\left( x,t\right) m\left( x,t\right) $, and it is
much easier to handle than the global interaction term. Unlike this, the
global interaction term is a highly unusual one in the theory of CIPs. This
is the integral operator 
\begin{equation}
\dint\limits_{\Omega }Y\left( x,y\right) m\left( y,t\right) dy  \label{1.1}
\end{equation}%
with the kernel $Y\left( x,y\right) \in L_{\infty }\left( \Omega \times
\Omega \right) $ \cite[formulas (3.1)]{Chow}. The presence of this term
causes a significant challenge here. We refer to subsection 2.3 and the
first item of Remark 3.2 (subsection 3.1) for further comments about (\ref%
{1.1}).

Since two parabolic equations forming the MFGS have two opposite directions
of time, then the classical theory of parabolic PDEs in inapplicable to the
MFGS. Indeed, until very recently uniqueness theorems for the MFGS were
unknown unless restrictive conditions are not imposed \cite{Bardi}. Recently
the author with coauthors has addressed the uniqueness question in four
works \cite{MFG1}-\cite{MFG4}, all of which are for the case of a single
measurement data. This was done via introducing either one extra
terminal/initial condition in the MFGS \cite{MFG1,MFG2,MFG4} or via
considering lateral Cauchy data and ignoring terminal and initial conditions 
\cite{MFG3}. Both Lipschitz and H\"{o}lder stability estimates for solutions
of corresponding problems for the MFGS were obtained in these references.
Such estimates for the MFGS were unknown prior\ \cite{MFG1}-\cite{MFG4}.

For the first time, the apparatus of Carleman estimates was introduced in
the MFG theory in \cite{MFG1}. This apparatus is also used in \cite{MFG2}-%
\cite{MFG4} as well as in this paper. While works \cite{MFG1}-\cite{MFG4}
consider the case when coefficients of the MFGS are known, in \cite{MFG5}
both H\"{o}lder and Lipschitz stability estimates are derived for some CIPs
for the MFGS with the final overdetermination. The latter means that the
initial and terminal conditions are given for both functions $u\left(
x,t\right) $ and $m\left( x,t\right) $,\ and also either partial or complete
lateral Cauchy data are given for these two functions. In the case of one
parabolic equation, a CIP, similar with the one of \cite{MFG5}, was
considered by the author in \cite{Klibt=T}. CIPs of \cite{MFG5} are
significantly different from the one considered in this paper.

Both the method of \cite{MFG5} and the technique of the current paper use a
significantly modified framework of the so-called Bukhgeim-Klibanov method,
which was first published in \cite{BukhKlib}. In this reference, the tool of
Carleman estimates was introduced in the field of Inverse Problems for the
first time, see, e.g. \cite{ImYam1,ImLY,Klib2006,Ksurvey,KL,YY,Yam} and
references cited therein for some follow up publications. The author with
coauthors has also extended the idea of \cite{BukhKlib} from the theory to
the so-called convexification globally convergent numerical method, see,
e.g. \cite{KLpar,KL} including the most recent works for problems of the
MFGS \cite{MFG7,MFG8}.

There is also a significant additional difficulty here, which we overcome.
This difficulty is due to the fact that the unknown coefficient of the MFGS
is involved in this system together with its first $x-$derivatives.
Nevertheless we work with the data resulting from a single measurement
event. We also refer to \cite{Klib2006} for a CIP for a hyperbolic equation
and to \cite{YY} for CIP for a parabolic equation. In both these cases
unknown coefficients are involved together with their first derivatives.
However, it is yet unclear whether techniques of \cite{Klib2006,YY} can be
adapted for our case.

The author is aware about only three previous analytical works of other
authors on CIPs for the MFGS \cite{ImLY,Liu1,Liu2}. Publications \cite%
{Liu1,Liu2} are the first theoretical results for CIPs for the MFGs. In \cite%
{Liu1,Liu2} uniqueness of the reconstruction of the interaction term $%
F\left( x,m\left( x,t\right) \right) $ in the HJB equation is proven in the
case when multiple initial conditions $m\left( x,0\right) $ are known. Using
the Carleman estimate of \cite{ImYam1,Yam} and the framework of \cite%
{BukhKlib}, Lipschitz stability estimate is proven in \cite{ImLY} for an
inverse source problem for the MFGS. It is well known that the question
about stability and uniqueness of a CIP can usually be reduced to the same
questions for an inverse source problem. However, although the problem
considered in \cite{ImLY} is similar with the one of this paper, our case of
the CIP\ for the MFGS, in which the unknown coefficient is involved together
with its first derivatives in the underlying PDEs as well as integral term (%
\ref{1.1}) is present, does not fit the framework of \cite{ImLY}, also, see
subsection 2.3. We refer to \cite{Chow,Ding,MFG7,MFG8} for numerical studies
of CIPs for the MFGS.

\textbf{Remark 1.1. }\emph{We are not concerned below with extra smoothness
conditions since traditionally they are not of a significant concern in the
field of CIPs, see, e.g. \cite{Nov}, \cite[Theorem 4.1]{Rom2}. }

All functions of this paper are real valued ones. We use below
Cauchy-Schwarz and Young's inequalities multiple times without further
mentioning. We formulate our CIP in section 2. In section 3, we formulate
our theorems. We prove these theorems in sections 4 and 5.

\section{Coefficient Inverse Problem}

\label{sec:2}

\subsection{Problem statement}

\label{sec:2.1}

Let $x=\left( x_{1},x_{2},...,x_{n}\right) $ denotes points in $\mathbb{R}%
^{n}$ and let $\overline{x}=\left( x_{2},...,x_{n}\right) .$ To make the
presentation simpler, we assume that our domain of interest $\Omega \subset 
\mathbb{R}^{n}$ is a rectangular prism. Let $a,b,B_{i}>0,i=1,...,n$ and $T>0$
be some numbers and $a<b$. We set:%
\begin{equation}
\left. 
\begin{array}{c}
\Omega =\left\{ x:a<x_{1}<b,-B_{i}<x_{i}<B_{i},i=2,...,n\right\} ,\text{ }
\\ 
\Omega _{1}=\left\{ \overline{x}:-B_{i}<x_{i}<B_{i},i=2,...,n\right\} , \\ 
\Gamma _{1}^{+}=\left\{ x\in \partial \Omega :x_{1}=b\right\} ,\Gamma
_{1}^{-}=\left\{ x\in \partial \Omega :x_{1}=a\right\} ,\Gamma _{1,T}^{\pm
}=\Gamma _{1}^{\pm }\times \left( 0,T\right) , \\ 
\Gamma _{i}^{\pm }=\left\{ x\in \partial \Omega :x_{i}=\pm B_{i}\right\} ,%
\text{ }\Gamma _{iT}^{\pm }=\Gamma _{i}^{\pm }\times \left( 0,T\right) ,%
\text{ }i=2,...,n, \\ 
Q_{T}=\Omega \times \left( 0,T\right) ,\text{ }S_{T}=\partial \Omega \times
\left( 0,T\right) .%
\end{array}%
\right.  \label{2.1}
\end{equation}%
Keeping in mind notations (\ref{2.1}),\emph{\ }denote 
\begin{equation*}
\left. 
\begin{array}{c}
H^{2,1}\left( Q_{T}\right) =\left\{ u:\left\Vert u\right\Vert
_{H^{2r,r}\left( Q_{T}\right) }^{2}=\dsum\limits_{\left\vert \alpha
\right\vert +2m\leq 2}\left\Vert D_{x}^{\alpha }\partial
_{t}^{m}u\right\Vert _{L_{2}\left( Q_{T}\right) }^{2}<\infty \right\} , \\ 
H^{2,1}\left( \Gamma _{iT}^{\pm }\right) =\left\{ 
\begin{array}{c}
u:\left\Vert u\right\Vert _{H^{2,1}\left( \Gamma _{iT}^{\pm }\right)
}^{2}=\dsum\limits_{j=1,j\neq i}^{n}\left\Vert u_{x_{j}}\right\Vert
_{L_{2}\left( \Gamma _{jT}^{\pm }\right) }^{2}+ \\ 
+\dsum\limits_{j,s=1,\left( j,s\right) \neq \left( i,i\right)
}^{n}\left\Vert u_{x_{j}x_{s}}\right\Vert _{L_{2}\left( \Gamma _{jT}^{\pm
}\right) }^{2}+\dsum\limits_{j=0}^{1}\left\Vert \partial
_{t}^{j}u\right\Vert _{L_{2}\left( \Gamma _{jT}^{\pm }\right) }^{2}<\infty ,%
\end{array}%
\right\} , \\ 
H^{2,1}\left( S_{T}\right) =\left\{ u:\left\Vert u\right\Vert
_{H^{2,1}\left( S_{T}\right) }^{2}=\dsum\limits_{i=1}^{n}\left\Vert
u\right\Vert _{H^{2,k}\left( \Gamma _{iT}^{\pm }\right) }^{2}<\infty
\right\} , \\ 
H^{1,0}\left( \Gamma _{iT}^{\pm }\right) =\left\{ u:\left\Vert u\right\Vert
_{H^{1,0}\left( \Gamma _{iT}^{\pm }\right) }^{2}=\dsum\limits_{j=1,j\neq
i}^{n}\left\Vert u_{x_{j}}\right\Vert _{L_{2}\left( \Gamma _{iT}^{\pm
}\right) }^{2}+\left\Vert u\right\Vert _{L_{2}\left( \Gamma _{iT}^{\pm
}\right) }^{2}<\infty \right\} , \\ 
i=1,...,n, \\ 
H^{1,0}\left( S_{T}\right) =\left\{ u:\left\Vert u\right\Vert
_{H^{1,0}\left( S_{T}\right) }^{2}=\dsum\limits_{i=1}^{n}\left\Vert
u\right\Vert _{H^{1,0}\left( \partial _{i}^{\pm }\Omega _{T}\right)
}^{2}<\infty \right\} .%
\end{array}%
\right.
\end{equation*}

We consider the MFGS of the second order in the following form \cite[%
formulas (3.1)]{Chow},

\begin{equation}
\left. 
\begin{array}{c}
u_{t}(x,t)+\Delta u(x,t){-k(x)(\nabla u(x,t))^{2}/2}+ \\ 
+\dint\limits_{\Omega }Y\left( x,y\right) m\left( y,t\right) dy+f\left(
x,t\right) m\left( x,t\right) =0,\text{ }\left( x,t\right) \in Q_{T}, \\ 
m_{t}(x,t)-\Delta m(x,t){-\func{div}(k(x)m(x,t)\nabla u(x,t))}=0,\text{ }%
\left( x,t\right) \in Q_{T},%
\end{array}%
\right.  \label{2.4}
\end{equation}%
where ${\nabla u=}\left( u_{x_{1}},...,u_{x_{n}}\right) ,$ the coefficient $%
k(x)\in C^{1}\left( \overline{\Omega }\right) $ and $f\left( x,t\right)
m\left( x,t\right) $ is the local interaction term. Thus, we use in the
second line of (\ref{2.4}) both global and local interaction terms. We
assume below that the function $f\left( x,t\right) $ is twice differentiable
with respect to $t$ and functions $f,f_{t},f_{tt}\in L_{\infty }\left(
Q_{T}\right) .$ The first equation (\ref{2.4}) is the
Hamilton-Jacobi-Bellman (HJB) equation and the second one is the
Fokker-Plank equation. The function $Y\left( x,y\right) \in L_{\infty
}\left( \Omega \times \Omega \right) $ is specified in subsections 2.4 and
2.5.

We now formulate two CIPs, which we study in this paper.

\textbf{Coefficient Inverse} \textbf{Problem }(CIP). \emph{Assume that
functions }$u,m\in C^{4}\left( \overline{Q}_{T}\right) $\emph{\ satisfy
equations (\ref{2.4}). Let }$t_{0}\in \left( 0,T\right) $\emph{\ be a
number. Let}%
\begin{equation}
\left. 
\begin{array}{c}
u\left( x,t_{0}\right) =u_{0}\left( x\right) ,\text{ }m\left( x,t_{0}\right)
=m_{0}\left( x\right) ,\text{ }x\in \Omega , \\ 
u\mid _{S_{T}}=g_{0}\left( x,t\right) ,\text{ }\partial _{n}u\mid
_{S_{T}}=g_{1}\left( x,t\right) , \\ 
m\mid _{S_{T}}=p_{0}\left( x,t\right) ,\text{ }\partial _{n}m\mid
_{S_{T}}=p_{1}\left( x,t\right) ,%
\end{array}%
\right.  \label{2.5}
\end{equation}%
\emph{where }$n\left( x\right) $\emph{\ is the unit outward looking normal
vector at the point }$x\in \partial \Omega .$\emph{\ Determine the
coefficient }$k\left( x\right) \in C^{1}\left( \overline{\Omega }\right) ,$%
\emph{\ assuming that the right hand sides of equalities (\ref{2.5}) are
known. }

\textbf{Remarks 2.1:}

\begin{enumerate}
\item \emph{The right hand sides of second and third lines of equalities (%
\ref{2.5}) are called \textquotedblleft lateral Cauchy data".}

\item \emph{Even though the lateral Cauchy data in (\ref{2.5}) are given on
the entire lateral boundary }$S_{T}$\emph{\ and represent, therefore,
complete lateral Cauchy data, we also consider below the case of incomplete
data, see item 1 of Remarks 3.3, in subsection 3.3. }
\end{enumerate}

\subsection{Some comments}

\label{sec:2.2}

This is a CIP\ with the single measurement data. In a real game scenario,
the input data (\ref{2.5}) can be obtained by conducting polls of game
players. In the case of the boundary data, polls are usually conducted not
just at the boundary but rather in a small neighborhood of the boundary.
Consider, for example the 1d case when $\Omega =\left( a,b\right) .$ Then
polling for both value function $u\left( x,t\right) $ and the density
function $m\left( x,t\right) $ is conducted for $x\in \left( a,a+\gamma
\right) \cup \left( b-\gamma ,b\right) ,$ where $\gamma >0$ is a small
number. Hence, it is possible to approximately figure out both Dirichlet and
Neumann boundary conditions at $x=a,b$ for both functions $u\left(
x,t\right) $ and $m\left( x,t\right) .$ To figure out functions $u\left(
x,t_{0}\right) ,m\left( x,t_{0}\right) ,$ one needs to arrange polling of
game players at a fixed moment of time $\left\{ t=t_{0}\right\} .$ Similar
considerations are valid for the $n-$D case, which we consider here. Our
input data (\ref{2.5}) require less information than the data of \cite[Model
1, Model 2]{Ding}, in which case the data are given for all $\left(
x,t\right) \in Q_{T}.$

In the case of one parabolic equation, uniqueness theorems for CIPs for the
situation when the data are known at a moment of time $t_{0}\in \left(
0,T\right) ,$ in addition to the incomplete lateral Cauchy data, where
proven by the method of \cite{BukhKlib} in a number of publications of the
author, see, e.g. \cite{Ksurvey}, \cite[section 3.4.2]{KL}. In \cite%
{ImYam1,Yam} Lipschitz stability estimates were proven for this problem.

\subsection{Comments about the global interaction term (\protect\ref{1.1})}

\label{sec:2.3}

It was pointed out in section 1 that this term is the one, which causes a 
\underline{significant challenge} of this paper. Indeed, it is the presence
of this term in MFGS (\ref{2.4}), which does not allow almost
\textquotedblleft automatic" applications of the previously developed theory
of CIPs for a single parabolic equation to CIPs for the MFGS.

In \cite{ImYam1} Lipschitz stability estimate was proven for a CIP for a
single parabolic equation, and that CIP is quite similar with ours, also,
see \cite{Yam}. In addition, the Lipschitz stability estimate for an analog
of our CIP (\ref{2.4})-(\ref{2.5}) for the MFGS is proven in the recent work 
\cite{ImLY} using a close analog of the idea of \cite{ImYam1}. Thus, it
seems to be, at least at the first glance, that the Lipschitz estimate can
indeed be proven for our CIP (\ref{2.4})-(\ref{2.5}) via basically the same
method as the one in \cite{ImYam1,ImLY,Yam}. However, the above mentioned
significant challenge of the MFGS, which is 
\begin{equation}
Y\left( x,y\right) \neq 0,  \label{200}
\end{equation}%
is not addressed in \cite{ImLY}. This indicates that, at least at the time
being, the Carleman Weight Function, which is used in \cite{ImYam1,ImLY,Yam}
to obtain the above mentioned Lipschitz stability estimates, is not capable
to handle the case (\ref{200}) for our CIP (\ref{2.4})-(\ref{2.5}).
Furthermore, the author is unaware about any ideas of obtaining Lipschitz
stability estimates for problem (\ref{2.4})-(\ref{2.5}) under condition (\ref%
{200}).

\textbf{Remark 2.2}. \emph{In fact, our H\"{o}lder stability estimates can
be arranged to become almost Lipschitz stability estimates, see item 3 in
Remarks 3.3 in subsection 3.3.}

The author cannot handle the case (\ref{200}) for a general function $%
Y\left( x,y\right) $. Nevertheless, the author can handle two important less
general forms of the function $Y\left( x,y\right) ,$ see subsections 2.4 and
2.5. The price we pay for being able to handle our two cases with condition (%
\ref{200}) is that we obtain H\"{o}lder stability estimates for our CIP
rather than the anticipated and more challenging Lipschitz stability. Still,
the author firmly believes that H\"{o}lder stability estimates of this paper
are better than the current absence of any stability estimates with
condition (\ref{200}), especially due to the numerical consideration in the
next two paragraphs. To work out these, we use a special Carleman Weight
Function, see subsection 3.1.

Finally, another important consideration in favor of the H\"{o}lder
stability estimates of this paper is the numerical one. Indeed, we have
already obtained numerical results via the convexification method for our
CIP (\ref{2.4})-(\ref{2.5}) \cite{MFG8}. These results have about the same
good quality as the ones in our previous works for a similar CIP for a
single parabolic equation \cite{KLpar}, \cite[chapter 9]{KL}. Note that in
our recent numerical work \cite{MFG7} a version of the convexification
method was applied to a different problem for the MFGS and good quality
results were obtained. The distinctive feature of the convexification is its
rigorously guaranteed global convergence, which is a rare case in the field
of CIPs.

In any version of the convexification method one basically adapts the proof
of the corresponding uniqueness/stability theorem to specific needs of the
convexification. This is exactly what is done in \cite{MFG8}. The presence
of the Carleman Weight Function in the numerical scheme is a crucial element
of the convexification. However, the Carleman Weight Function of \cite%
{ImYam1,ImLY,Yam} depends on two large parameters, which means that it
changes too rapidly. The latter is inconvenient for the numerical
implementation. On the other hand, in \cite{MFG8} the Carleman Weight
Function of subsection 3.1 is used, that function depends on a single large
parameter, and our experience of \cite{KLpar}, \cite[chapter 9]{KL}
demonstrates that this function is truly convenient to work with in
numerical studies.

\subsection{The first form of the kernel $Y\left( x,y\right) $ of the
integral operator in (\protect\ref{2.4})}

\label{sec:2.4}

It was pointed out in \cite[page 2653]{LiuOsher} that a popular choice of
the kernel $Y\left( x,y\right) $ of the integral operator in the first
equation (\ref{2.4}) is the product of $n$ Gaussians,%
\begin{equation}
\dprod\limits_{i=1}^{n}\exp \left( -\frac{\left( x_{i}-y_{i}\right) ^{2}}{%
2\sigma _{i}^{2}}\right)  \label{2.0}
\end{equation}%
with some numbers $\sigma _{i}>0.$ Since Gaussians approximate the $\delta -$%
function in the distribution sense, then this justifies our choice of the
function $Y\left( x,y\right) $ in (\ref{1.1}) in the form, which is somewhat
more general than (\ref{2.0}), i.e. 
\begin{equation}
\left. 
\begin{array}{c}
Y\left( x,y\right) =\delta \left( x_{1}-y_{1}\right) \overline{Y}\left( 
\overline{x},\overline{y}\right) , \\ 
\text{ }\overline{Y}\in L_{\infty }\left( \Omega _{1}\times \Omega
_{1}\right) ,\text{ }\left\Vert \text{ }\overline{Y}\right\Vert _{L_{\infty
}\left( \Omega _{1}\times \Omega _{1}\right) }\leq N_{1}.%
\end{array}%
\right.  \label{2.2}
\end{equation}%
In (\ref{2.2}) $\delta \left( x_{1}-y_{1}\right) $ is the $\delta -$function
and $N_{1}>0$ is a number. By (\ref{2.2}) the integral in the first equation
(\ref{2.4}) becomes%
\begin{equation}
\dint\limits_{\Omega }Y\left( x,y\right) m\left( y,t\right)
dy=\dint\limits_{\Omega _{1}}\overline{Y}\left( \overline{x},\overline{y}%
\right) m\left( x_{1},\overline{y},t\right) d\overline{y},\text{ }x=\left(
x_{1},\overline{x}\right) \in \Omega ,t\in \left( 0,T\right) .  \label{2.3}
\end{equation}%
Thus, the integral operator in the right hand side of (\ref{2.3}) represents
local interaction with respect to $x_{1}$ and global interaction with
respect to $x_{2},...,x_{n}.$

\subsection{The second form of the kernel $Y\left( x,y\right) $ of the
integral operator in (\protect\ref{2.4})}

\label{sec:2.5}

For $z\in \mathbb{R}$ let $H\left( z\right) $ be the Heaviside function,%
\begin{equation*}
H\left( z\right) =\left\{ 
\begin{array}{c}
1,z>0, \\ 
0,z<0.%
\end{array}%
\right.
\end{equation*}%
Consider the function $Y\left( x,y\right) $ in the form%
\begin{equation}
\left. 
\begin{array}{c}
Y\left( x,y\right) =H\left( y_{1}-x_{1}\right) \overline{Y}\left( x,y\right)
,\text{ }\left( x,y\right) \in \Omega \times \Omega , \\ 
\text{ }\overline{Y}\in L_{\infty }\left( \Omega \times \Omega \right) ,%
\text{ }\left\Vert \text{ }\overline{Y}\right\Vert _{L_{\infty }\left(
\Omega \times \Omega \right) }\leq N_{1}.%
\end{array}%
\right.  \label{2.6}
\end{equation}%
Then by (\ref{2.1}) and (\ref{2.6}) the integral in the first equation (\ref%
{2.4}) becomes%
\begin{equation}
\left. 
\begin{array}{c}
\dint\limits_{\Omega }Y\left( x,y\right) m\left( y,t\right) dy= \\ 
=\dint\limits_{\Omega _{1}}\left( \dint\limits_{x_{1}}^{b}\overline{Y}\left(
x_{1},\overline{x},y_{1},\overline{y}\right) m\left( y_{1},\overline{y}%
,t\right) dy_{1}\right) d\overline{y},\text{ }x=\left( x_{1},\overline{x}%
\right) \in \Omega .%
\end{array}%
\right.  \label{2.7}
\end{equation}%
Therefore, the integral operator in (\ref{2.7}) represents the global
interaction with respect to all $n$ spatial variables $x_{1},...,x_{n}.$

\section{Formulations of Theorems}

\label{sec:3}

\textbf{Remark 3.1.} \emph{Let }$t_{0}\in \left( 0,T\right) $\emph{\ be the
number in (\ref{2.5}). Without any loss of generality we set everywhere
below }$t_{0}=T/2.$\emph{\ }

\subsection{Two Carleman estimates}

\label{sec:3 .1}

Let $\alpha >0$ and $\lambda \geq 1$ be two parameters, which we will choose
later, and $\lambda $ will be sufficiently large. The Carleman Weight
Function in our case is $\varphi _{\lambda }\left( x_{1},t\right) ,$%
\begin{equation}
\varphi _{\lambda }\left( x_{1},t\right) =\exp \left[ 2\lambda \left(
x_{1}^{2}-\alpha \left( t-T/2\right) ^{2}\right) \right] .  \label{3.01}
\end{equation}%
Let the number $\varepsilon \in \left( 0,T/2\right) .$ Denote%
\begin{equation}
Q_{\varepsilon ,T}=\Omega \times \left( \varepsilon ,T-\varepsilon \right) .
\label{3.001}
\end{equation}%
It follows from (\ref{2.1}), (\ref{3.01}) and (\ref{3.001}) that%
\begin{equation}
\left. 
\begin{array}{c}
\max_{\overline{Q}_{T}}\varphi _{\lambda }\left( x_{1},t\right) =\varphi
_{\lambda }\left( b,T/2\right) =e^{2\lambda b^{2}}, \\ 
\min_{\overline{Q}_{T,\varepsilon }}\varphi _{\lambda }\left( x_{1},t\right)
=\exp \left[ 2\lambda \left( a^{2}-\alpha \left( T/2-\varepsilon \right)
^{2}\right) \right] .%
\end{array}%
\right.  \label{3.2}
\end{equation}

\textbf{Remarks 3.2:}

\begin{enumerate}
\item \textbf{\ }\emph{The Carleman Weight Function (\ref{3.01}) for
parabolic operators }$\partial _{t}\pm \Delta $ \emph{depends only on one
spatial variable }$x_{1}$\emph{\ and has only one large parameter }$\lambda $
\emph{rather than two large parameters of} \emph{\cite%
{ImYam1,ImLY,Ksurvey,YY,Yam}, \cite[section 2.3]{KL}, \cite[Chapter 4, \S 1]%
{LRS} and other references.\ This Carleman Weight Function was previously
used in the works of the author with coauthors \cite{KLpar,MFG3,MFG8}, \cite[%
formula (9.20)]{KL}.\ In particular, publications \cite{KLpar}, \cite[%
chapter 9]{KL} are about numerical studies of a similar CIP for one
parabolic equation via the convexification method, and the most recent work 
\cite{MFG8} is about numerical studies of our CIP (\ref{2.4})-(\ref{2.5})
using the convexification method. As stated in subsection 2.3, function (\ref%
{3.01}) is more convenient for the numerical implementation in the
convexification method than those depending on two large parameters. The
author is unaware about other publications, in which Carleman Weight
Function (\ref{3.01}) is used. }

\item \emph{The reason why we choose (\ref{3.01}) is that it allows us to
work with the integral term in (\ref{2.4}) for cases (\ref{2.3}), (\ref{2.7}%
). }

\item \emph{The reason of our assumption that the domain }$\Omega $\emph{\
in (\ref{2.1}) is a rectangular prism is that Theorem 3.1 is proven in \cite%
{MFG3} only for this case. It is not yet clear whether this theorem is valid
for a more general domain }$\Omega .$ \emph{Still, it is clear from (\ref%
{2.3}) and (\ref{2.7}) that }$\Omega $\emph{\ should be a cylindrical domain.%
}

\item \emph{In addition to item 3, based on his experience, the author
believes that a non very general shape of the underlying domain is rarely a
primary concern when working with CIPs, since CIPs are very difficult ones
in their own rights anyway, see, e.g. \cite[sections 2.3, 2.4]{Ksurvey}, 
\cite{KLpar}, \cite[sections 3.2, chapters 7-10]{KL}, \cite[\S 1 in Chapter 4%
]{LRS}.}
\end{enumerate}

\textbf{Theorem 3.1 }\cite{MFG3}\textbf{.}\emph{\ Let conditions (\ref{2.1})
and (\ref{3.01}) hold. Then there exists a sufficiently large number }$%
\lambda _{0}=\lambda _{0}\left( \alpha ,\Omega ,T\right) \geq 1$\emph{\ and
a number }$C_{0}=C_{0}\left( \alpha ,\Omega ,T\right) >0,$\emph{\ both
numbers depending only on listed parameters, such that the following two
Carleman estimates hold: }%
\begin{equation}
\left. 
\begin{array}{c}
\dint\limits_{Q_{T}}\left( u_{t}\pm \Delta u\right) ^{2}\varphi _{\lambda
}dxdt\geq \left( C_{0}/\lambda \right) \dint\limits_{Q_{T}}\left(
u_{t}^{2}+\dsum\limits_{i,j=1}^{n}u_{x_{i}x_{j}}^{2}\right) \varphi
_{\lambda }dxdt+ \\ 
+C_{0}\dint\limits_{Q_{T}}\left( \lambda \left( \nabla u\right) ^{2}+\lambda
^{3}u^{2}\right) \varphi _{\lambda }dxdt- \\ 
-C_{0}\left( \left\Vert \partial _{n}u\right\Vert _{H^{1,0}\left(
S_{T}\right) }^{2}+\left\Vert u\right\Vert _{H^{2,1}\left( S_{T}\right)
}^{2}\right) e^{3\lambda b^{2}}- \\ 
-C_{0}\left( \left\Vert u\left( x,T\right) \right\Vert _{H^{1}\left( \Omega
\right) }^{2}+\left\Vert u\left( x,0\right) \right\Vert _{H^{1}\left( \Omega
\right) }^{2}\right) \exp \left( -2\lambda \left( \alpha
T^{2}/4-b^{2}\right) \right) ,\text{ } \\ 
\forall \lambda \geq \lambda _{0},\text{ }\forall u\in H^{3}\left(
Q_{T}\right) .%
\end{array}%
\right.  \label{3.04}
\end{equation}%
\emph{Suppose now that }%
\begin{equation}
u\mid _{S_{T}\diagdown \left( \Gamma _{1T}^{+}{}\right) }=0.  \label{3.05}
\end{equation}%
\emph{Then the third line of (\ref{3.04}) becomes}%
\begin{equation}
-C\left( \left\Vert \partial _{n}u\right\Vert _{H^{1,0}\left( \Gamma
_{1T}^{+}\right) }^{2}+\left\Vert u\right\Vert _{H^{2,1}\left( \Gamma
_{1T}^{+}\right) }^{2}\right) e^{3\lambda b^{2}}.  \label{3.06}
\end{equation}

The implication of (\ref{3.06}) from (\ref{3.05}) can be obtained by a
slight modification of the proof of \cite{MFG3}.

\subsection{Estimates of some integrals}

\label{sec:3.2}

In this subsection we estimate from the above weighted $L_{2}\left(
Q_{T}\right) $ norms of terms of terms in the right hand sides of (\ref{2.3}%
) and (\ref{2.6}) as well as one more integral. The weight is the Carleman
Weight Function (\ref{3.01}). In this subsection $\widetilde{C}=\widetilde{C}%
\left( N_{1},Q_{T}\right) >0$ denotes different numbers depending only on
listed parameters.

\textbf{Lemma 3.1}. \ \emph{Let }$\overline{Y}\left( \overline{x},\overline{y%
}\right) $\emph{\ be the function in (\ref{2.2}), (\ref{2.3}). Then the
following estimate is valid}%
\begin{equation}
\left. 
\begin{array}{c}
\dint\limits_{Q_{T}}\left( \dint\limits_{\Omega _{1}}\overline{Y}\left( 
\overline{x},\overline{y}\right) h\left( x_{1},\overline{y},t\right) d%
\overline{y}\right) ^{2}\varphi _{\lambda }\left( x,t\right) dxdt\leq \\ 
\leq \widetilde{C}\dint\limits_{Q_{T}}h^{2}\varphi _{\lambda }dxdt,\text{ }%
\forall h\in L_{2}\left( Q_{T}\right) ,\text{ }\forall \lambda >0.%
\end{array}%
\right.  \label{3.7}
\end{equation}

The proof of this Lemma is obvious since the function $\varphi _{\lambda
}\left( x_{1},t\right) $ is independent on $\overline{x}.$

\bigskip \textbf{Lemma 3.2}. \ \emph{Let }$Y\left( x,y\right) $\emph{\ be
the function in (\ref{2.6}). Then the following estimate is valid}%
\begin{equation}
\left. 
\begin{array}{c}
\dint\limits_{Q_{T}}\left( \dint\limits_{\Omega }Y\left( x,y\right) h\left(
y,t\right) dy\right) ^{2}\varphi _{\lambda }\left( x,t\right) dxdt\leq \\ 
\leq \widetilde{C}\dint\limits_{Q_{T}}h^{2}\varphi _{\lambda }dxdt,\text{ }%
\forall h\in L_{2}\left( Q_{T}\right) ,\text{ }\forall \lambda >0.%
\end{array}%
\right.  \label{3.8}
\end{equation}

\textbf{Proof}. By (\ref{2.7}), Cauchy-Schwarz inequality and Fubini theorem%
\begin{equation}
\left. 
\begin{array}{c}
\dint\limits_{Q_{T}}\left( \dint\limits_{\Omega }Y\left( x,y\right) h\left(
y,t\right) dy\right) ^{2}\varphi _{\lambda }\left( x,t\right) dxdt\leq \\ 
\leq \widetilde{C}\dint\limits_{Q_{T}}\left[ \dint\limits_{\Omega
_{1}}\left( \dint\limits_{x_{1}}^{b}h^{2}\left( y_{1},\overline{y},t\right)
dy_{1}\right) d\overline{y}\right] \varphi _{\lambda }\left( x_{1},t\right)
dxdt= \\ 
=\widetilde{C}\dint\limits_{0}^{T}\dint\limits_{\Omega _{1}}\left\{
\dint\limits_{a}^{b}\left[ \dint\limits_{\Omega _{1}}\left(
\dint\limits_{x_{1}}^{b}h^{2}\left( y_{1},\overline{y},t\right)
dy_{1}\right) d\overline{y}\right] \varphi _{\lambda }\left( x_{1},t\right)
dx_{1}\right\} d\overline{x}dt\leq \\ 
\leq \widetilde{C}\dint\limits_{0}^{T}\dint\limits_{\Omega _{1}}\left[
\dint\limits_{a}^{b}\left( \dint\limits_{x_{1}}^{b}h^{2}\left( y_{1},%
\overline{y},t\right) dy_{1}\right) \varphi _{\lambda }\left( x_{1},t\right)
dx_{1}\right] d\overline{y}dt.%
\end{array}%
\right.  \label{3.9}
\end{equation}%
We obviously have 
\begin{equation*}
\dint\limits_{a}^{b}\left( \dint\limits_{x_{1}}^{b}f\left(
x_{1},y_{1}\right) dy_{1}\right) dx_{1}=\dint\limits_{a}^{b}\left(
\dint\limits_{a}^{y_{1}}f\left( x_{1},y_{1}\right) dx_{1}\right) dy_{1},%
\text{ }\forall f\in L_{1}\left( \left( a,b\right) \times \left( a,b\right)
\right) .
\end{equation*}%
Hence, using the fact that the $\varphi _{\lambda }\left( x_{1},t\right) $
is increasing with respect to $x_{1}\in \left( a,b\right) ,$ we obtain in
the last integral of (\ref{3.9})%
\begin{equation*}
\left. 
\begin{array}{c}
\dint\limits_{a}^{b}\left( \dint\limits_{x_{1}}^{b}h^{2}\left( y_{1},%
\overline{y},t\right) dy_{1}\right) \varphi _{\lambda }\left( x_{1},t\right)
dx_{1}=\dint\limits_{a}^{b}\left( \dint\limits_{a}^{y_{1}}\varphi _{\lambda
}\left( x_{1},t\right) dx_{1}\right) h^{2}\left( y_{1},\overline{y},t\right)
dy_{1}\leq \\ 
\leq \dint\limits_{a}^{b}\left( \dint\limits_{a}^{y_{1}}\varphi _{\lambda
}\left( y_{1},t\right) dx_{1}\right) h^{2}\left( y_{1},\overline{y},t\right)
dy_{1}\leq \left( b-a\right) \dint\limits_{a}^{b}h^{2}\left( y_{1},\overline{%
y},t\right) \varphi _{\lambda }\left( y_{1},t\right) dy_{1}.%
\end{array}%
\right.
\end{equation*}%
Hence,%
\begin{equation*}
\left. 
\begin{array}{c}
\dint\limits_{0}^{T}\dint\limits_{\Omega _{1}}\left[ \dint\limits_{a}^{b}%
\left( \dint\limits_{x_{1}}^{b}h^{2}\left( y_{1},\overline{y},t\right)
dy_{1}\right) \varphi _{\lambda }\left( x_{1},t\right) dx_{1}\right] d%
\overline{y}dt\leq \\ 
\leq \left( b-a\right) \dint\limits_{0}^{T}\left( \dint\limits_{\Omega
}h^{2}\left( x,t\right) \right) \varphi _{\lambda }\left( x,t\right)
dxdt=\left( b-a\right) \dint\limits_{Q_{T}}h^{2}\varphi _{\lambda }dxdt.%
\end{array}%
\right.
\end{equation*}%
Comparing this with (\ref{3.9}), we obtain the target estimate (\ref{3.8}).
\ $\square $

\textbf{Lemma 3.3 }(\cite{Ksurvey}, \textbf{\ }\cite[Lemma 3.1.1]{KL})%
\textbf{. }\emph{The following inequality holds}%
\begin{equation}
\left. 
\begin{array}{c}
\dint\limits_{Q_{T}}\left( \dint\limits_{T/2}^{t}h\left( x,\tau \right)
d\tau \right) ^{2}\varphi _{\lambda }dxdt\leq \\ 
\leq C_{1}\left( 1/\lambda \right) \dint\limits_{Q_{T}}\left( h^{2}\varphi
_{\lambda }\right) \left( x,t\right) dxdt,\text{ }\forall \lambda >0,\text{ }%
\forall h\in L_{2}\left( Q_{T}\right) ,%
\end{array}%
\right.  \label{3.10}
\end{equation}%
\emph{\ where the number }$C_{1}=C_{1}\left( \alpha \right) >0$\emph{\
depends only on the parameter }$\alpha $\emph{\ in (\ref{3.01}).}

\subsection{H\"{o}lder stability estimates and uniqueness}

\label{sec:3.3}

Let $N_{2},N_{3},N_{4}>0$ be three numbers. Recall that the number $N_{1}>0$
was introduced in (\ref{2.2}). Denote%
\begin{equation}
\left. 
\begin{array}{c}
S_{1}\left( N_{2}\right) =\left\{ v\in C^{4}\left( \overline{Q}_{T}\right)
:\left\Vert v\right\Vert _{C^{4}\left( \overline{Q}_{T}\right) }\leq
N_{2}\right\} , \\ 
S_{2}\left( N_{3}\right) =\left\{ k\in C^{1}\left( \overline{\Omega }\right)
:\left\Vert k\right\Vert _{C^{1}\left( \overline{\Omega }\right) }\leq
N_{3}\right\} , \\ 
N=\max \left( N_{1},N_{2},N_{3}\right) .%
\end{array}%
\right.  \label{4.1}
\end{equation}%
Suppose that we have two triples 
\begin{equation}
\left( u_{i},m_{i},k_{i}\right) \in S_{1}^{2}\left( N_{2}\right) \times
S_{2}\left( N_{3}\right) ,i=1,2  \label{4.2}
\end{equation}%
satisfying the following analogs of conditions (\ref{2.5}), see Remark 3.1:%
\begin{equation}
\left. 
\begin{array}{c}
u_{i}\left( x,T/2\right) =u_{0,i}\left( x\right) ,\text{ }m_{i}\left(
x,T/2\right) =m_{0,i}\left( x\right) ,\text{ }x\in \Omega ,\text{ }i=1,2, \\ 
u_{i}\mid _{S_{T}}=g_{0,i}\left( x,t\right) ,\text{ }\partial _{n}u_{i}\mid
_{S_{T}}=g_{1,i}\left( x,t\right) ,\text{ }i=1,2, \\ 
m_{i}\mid _{S_{T}}=p_{0}\left( x,t\right) ,\text{ }\partial _{n}m_{i}\mid
_{S_{T}}=p_{1,i}\left( x,t\right) ,\text{ }i=1,2.%
\end{array}%
\right.  \label{4.3}
\end{equation}%
Denote%
\begin{equation}
\left. 
\begin{array}{c}
\widetilde{u}=u_{1}-u_{2},\text{ }\widetilde{m}=m_{1}-m_{2},\text{ }%
\widetilde{k}=k_{1}-k_{2}, \\ 
\widetilde{u}_{0}=u_{0,1}-u_{0,2},\text{ }\widetilde{m}_{0}=m_{0,1}-m_{0,2},
\\ 
\widetilde{g}_{0}=g_{0,1}-g_{0,2},\text{ }\widetilde{g}_{1}=g_{1,1}-g_{1,2},%
\text{ }\widetilde{p}_{0}=p_{0,1}-p_{0,2},\text{ }\widetilde{p}%
_{1}=p_{1,1}-p_{1,2}.%
\end{array}%
\right.  \label{4.4}
\end{equation}%
Using (\ref{4.1})-(\ref{4.4}) and triangle inequality, we obtain%
\begin{equation}
\left. 
\begin{array}{c}
\left( \widetilde{u},\widetilde{m},\widetilde{k}\right) \in S_{1}^{2}\left(
2N_{2}\right) \times S_{2}\left( 2N_{3}\right) , \\ 
\widetilde{u}\left( x,T/2\right) =\widetilde{u}_{0}\left( x\right) ,\text{ }%
\widetilde{m}\left( x,T/2\right) =\widetilde{m}_{0}\left( x\right) ,\text{ }%
x\in \Omega , \\ 
\widetilde{u}\mid _{S_{T}}=\widetilde{g}_{0}\left( x,t\right) ,\text{ }%
\partial _{n}\widetilde{u}\mid _{S_{T}\diagdown \Gamma _{1T}^{-}}=\widetilde{%
g}_{1}\left( x,t\right) , \\ 
\widetilde{m}\mid _{S_{T}}=\widetilde{p}_{0}\left( x,t\right) ,\text{ }%
\partial _{n}\widetilde{m}\mid _{S_{T}{}^{-}}=\widetilde{p}_{1}\left(
x,t\right) .%
\end{array}%
\right.  \label{4.5}
\end{equation}

\textbf{Theorem 3.2.} \emph{Let either conditions (\ref{2.2}), (\ref{2.3})
or conditions (\ref{2.6}), (\ref{2.7}) be in place. Suppose that there exist
two triples }$\left( u_{i},m_{i},k_{i}\right) ,$\emph{\ }$i=1,2$\emph{\
satisfying equations (\ref{2.4}) and conditions (\ref{4.1})-(\ref{4.5}).
Assume that }%
\begin{equation}
\frac{1}{2}\left\vert \nabla u_{0,1}\left( x\right) \right\vert ^{2}\geq c,%
\text{ }x\in \overline{\Omega },  \label{4.6}
\end{equation}%
\emph{where }$c>0$\emph{\ is a number. Also, let the function }$f\left(
x,t\right) $\emph{\ in the first equation (\ref{2.4}) be such that there
exist its }$t-$\emph{derivatives }$f,f_{t},f_{tt}\in L_{\infty }\left(
Q_{T}\right) $\emph{\ and }$\left\Vert \partial _{t}^{k}f\right\Vert
_{L_{\infty }\left( Q_{T}\right) }\leq N,$ $k=0,1,2.$\emph{\ Let }$\delta >0$%
\emph{\ be a sufficiently small number. Assume that} 
\begin{equation}
\left. 
\begin{array}{c}
\left\Vert \widetilde{u}_{0}\right\Vert _{H^{1}\left( \Omega \right)
},\left\Vert \widetilde{m}_{0}\right\Vert _{H^{1}\left( \Omega \right) }\leq
\delta , \\ 
\left\Vert \partial _{t}^{s}\widetilde{g}_{0}\right\Vert _{H^{2,1}\left(
S_{T}{}^{-}\right) },\left\Vert \partial _{t}^{s}\widetilde{p}%
_{0}\right\Vert _{H^{2,1}\left( S_{T}\right) }\leq \delta ,\text{ }s=0,1,2,
\\ 
\left\Vert \partial _{t}^{s}\widetilde{g}_{1}\right\Vert _{H^{1,0}\left(
S_{T}\right) },\left\Vert \partial _{t}^{s}\widetilde{p}_{1}\right\Vert
_{H^{1,0}\left( S_{T}\right) }\leq \delta ,\text{ }s=0,1,2.%
\end{array}%
\right.  \label{4.7}
\end{equation}%
\emph{Let }$\rho \in \left( 0,1\right) $ \emph{be an arbitrary number. Then
for every }$\varepsilon $ \emph{satisfying}%
\begin{equation}
\frac{T}{2}\left( 1-\sqrt{\rho }\right) <\varepsilon <\frac{T}{2}.
\label{4.69}
\end{equation}%
\emph{\ there exists a sufficiently small number }%
\begin{equation}
\delta _{0}=\delta _{0}\left( N,\varepsilon ,\Omega ,T,c,\rho \right) \in
\left( 0,1\right) ,\text{ }  \label{4.70}
\end{equation}%
\emph{\ and a number }$C=C\left( N,\varepsilon ,\Omega ,T,c\right) >0$\emph{%
, both\ numbers depending only on listed parameters,\ such that} 
\begin{equation}
\left\Vert \partial _{t}^{s}\widetilde{u}\right\Vert _{H^{2,1}\left(
Q_{\varepsilon ,T}\right) },\left\Vert \partial _{t}^{s}\widetilde{m}%
\right\Vert _{H^{2,1}\left( Q_{\varepsilon ,T}\right) },\left\Vert 
\widetilde{k}\right\Vert _{L_{2}\left( \Omega \right) }\leq C\delta ^{1-\rho
},\text{ }\forall \delta \in \left( 0,\delta _{0}\right) ,\text{ }s=0,1,2.
\label{4.8}
\end{equation}%
\emph{Also, CIP (\ref{2.4})-(\ref{2.5}) has at most one solution }$\left(
u,m,k\right) \in S_{1}^{2}\left( N_{2}\right) \times S_{2}\left(
N_{3}\right) .$

\textbf{Theorem 3.3.} \emph{Assume that all conditions of Theorem 3.2 hold,
except that the normal derivatives }$\partial _{n}u_{i}$\emph{\ and }$%
\partial _{n}m_{i},$\emph{\ }$i=1,2$\emph{\ are unknown at }$S_{T}\diagdown
\Gamma _{1T}^{+}.$\emph{\ On the other hand, assume that (\ref{4.7}) is
replaced with\ }%
\begin{equation}
\left. 
\begin{array}{c}
\left\Vert \widetilde{u}_{0}\right\Vert _{H^{2}\left( \Omega \right)
},\left\Vert \widetilde{m}_{0}\right\Vert _{H^{1}\left( \Omega \right) }\leq
\delta , \\ 
\widetilde{g}_{0}\left( x,t\right) =\widetilde{p}_{0}\left( x,t\right)
=0,\left( x,t\right) \in S_{T}\diagdown \Gamma _{1,T}^{+}, \\ 
\left\Vert \partial _{t}^{s}\widetilde{g}_{1}\right\Vert _{H^{1,0}\left(
\Gamma _{1,T}^{+}\right) },\left\Vert \partial _{t}^{s}\widetilde{p}%
_{1}\right\Vert _{H^{1,0}\left( \Gamma _{1,T}^{+}\right) }\leq \delta ,\text{
}s=0,1,2.%
\end{array}%
\right.  \label{4.80}
\end{equation}%
\emph{Then estimate (\ref{4.8}) still holds as well as uniqueness. }

\textbf{Remarks 3.3}:

\begin{enumerate}
\item \emph{Therefore, it follows from (\ref{4.3}), (\ref{4.4}), (\ref{4.80}%
) and Theorem 3.3 that if Dirichlet boundary conditions for two pairs }$%
\left( u_{1},m_{1}\right) $\emph{\ and }$\left( u_{2},m_{2}\right) $\emph{\
coincide at }$S_{T}\diagdown \Gamma _{T}^{+},$\emph{\ then stability and
uniqueness for our CIP still hold, even if the Neumann boundary conditions
at }$S_{T}\diagdown \Gamma _{T}^{+}$\emph{\ are unknown for these two pairs.
In other words, Theorem 3.3 works for the case of incomplete lateral Cauchy
data for our CIP.}

\item \emph{Everywhere below }$C=C\left( N,\varepsilon ,\Omega ,T,c\right)
>0 $\emph{\ denotes different numbers depending only on listed parameters.
Also, in the course of the proof of Theorem 3.2 we will choose the parameter 
}$\alpha >0$ \emph{in the Carleman Weight Function (\ref{3.01}), so that it
will depend only on }$\Omega $, $T$\emph{\ and }$\varepsilon ,$\emph{\ i.e. }%
$\alpha =\alpha \left( \Omega ,T,\varepsilon \right) >0.$\emph{\ Hence, we
do not indicate the dependence of some parameters on }$\alpha $ \emph{in the
proof of Theorem 3.2.}

\item \emph{Note that if }$\rho \approx 0,$ \emph{then (\ref{4.8}) is almost
Lipschitz stability estimate, in which case one must have }$C\delta $\emph{\
instead of }$C\delta ^{1-\rho }$\emph{\ in (\ref{4.8}).}

\item A \emph{condition like (\ref{4.6}) is always imposed in publications
about the Bukhgeim-Klibanov method, including the above cited ones \cite%
{BukhKlib,ImYam1,ImLY,Klib2006,Ksurvey,KL,YY,Yam}.}
\end{enumerate}

\section{Proof of Theorem 3.2}

\label{sec:4}

In this proof the integral term in the first equation (\ref{2.4}) is:%
\begin{equation}
\left. 
\begin{array}{c}
\dint\limits_{\Omega }Y\left( x,y\right) m\left( y,t\right) dy= \\ 
=\left\{ 
\begin{array}{c}
\dint\limits_{\Omega _{1}}\overline{Y}\left( \overline{x},\overline{y}%
\right) m\left( x_{1},\overline{y},t\right) d\overline{y},\text{ in the case
(\ref{2.3}),} \\ 
\dint\limits_{\Omega _{1}}\left( \dint\limits_{x_{1}}^{b}\overline{Y}\left(
x_{1},\overline{x},y_{1},\overline{y}\right) m\left( y_{1},\overline{y}%
,t\right) dy_{1}\right) d\overline{y}\text{ in the case (\ref{2.7}).}%
\end{array}%
\right.%
\end{array}%
\right.  \label{4.800}
\end{equation}

Subtract equations (\ref{2.4}) for the triple $\left(
u_{2},m_{2},k_{2}\right) $ from the same equations for the triple $\left(
u_{1},m_{1},k_{1}\right) .$ Use the formula 
\begin{equation*}
y_{1}z_{1}-y_{2}z_{2}=\widetilde{y}z_{1}+y_{2}\widetilde{z},\text{ }\forall
y_{1},y_{2},z_{1},z_{2}\in \mathbb{R},
\end{equation*}%
where $\widetilde{y}=y_{1}-y_{2}$ and $\widetilde{z}=z_{1}-z_{2}.$ Using
notations (\ref{4.4}), we obtain%
\begin{equation}
\left. 
\begin{array}{c}
\widetilde{u}_{t}+\Delta \widetilde{u}+\dint\limits_{\Omega }Y\left(
x,y\right) \widetilde{m}\left( y,t\right) dy+f\left( x,t\right) \widetilde{m}%
\left( x,t\right) - \\ 
-k_{2}\left( x\right) \nabla \widetilde{u}\left( \nabla u_{1}+\nabla
u_{2}\right) /2=\widetilde{k}\left( x\right) \left( \nabla u_{1}\right)
^{2}/2,\text{ }\left( x,t\right) \in Q_{T},%
\end{array}%
\right.  \label{4.9}
\end{equation}%
\begin{equation}
\left. 
\begin{array}{c}
\widetilde{m}_{t}-\Delta \widetilde{m}-\func{div}\left( k_{2}\left( x\right) 
\widetilde{m}\nabla u_{1}\right) -\func{div}\left( k_{2}\left( x\right)
m_{2}\nabla \widetilde{u}\right) = \\ 
=\func{div}\left( \widetilde{k}\left( x\right) m_{1}\nabla u_{1}\right) ,%
\text{ }\left( x,t\right) \in Q_{T}.%
\end{array}%
\right.  \label{4.90}
\end{equation}
Denote 
\begin{equation}
v\left( x,t\right) =\widetilde{u}_{t}\left( x,t\right) ,\text{ }q\left(
x,t\right) =\widetilde{m}_{t}\left( x,t\right) .  \label{4.10}
\end{equation}%
Then by (\ref{4.5}) and (\ref{4.10})%
\begin{equation}
\widetilde{u}\left( x,t\right) =\dint\limits_{T/2}^{t}v\left( x,\tau \right)
d\tau +\widetilde{u}_{0}\left( x\right) ,\text{ }\widetilde{m}\left(
x,t\right) =\dint\limits_{T/2}^{t}q\left( x,\tau \right) d\tau +\widetilde{m}%
_{0}\left( x\right) .  \label{4.11}
\end{equation}%
Also, setting in (\ref{4.9}) $t=T/2$ and using (\ref{4.3}), (\ref{4.5}) and (%
\ref{4.10}), we obtain 
\begin{equation}
\widetilde{k}\left( x\right) =2\left( \nabla u_{0,1}\left( x\right) \right)
^{-2}v\left( x,T/2\right) +F\left( x\right) ,  \label{4.12}
\end{equation}%
where 
\begin{equation}
\left. 
\begin{array}{c}
F\left( x\right) = \\ 
=2\left( \nabla u_{0,1}\left( x\right) \right) ^{-2}\left( \Delta \widetilde{%
u}_{0}+\dint\limits_{\Omega }Y\left( x,y\right) \widetilde{m}_{0}\left(
y\right) dy+f\left( x,0\right) \widetilde{m}_{0}\left( x\right) \right) - \\ 
-\left( \nabla u_{0,1}\left( x\right) \right) ^{-2}k_{2}\left( x\right)
\nabla \widetilde{u}_{0}\nabla \left( u_{0,1}+u_{0,2}\right) \left( x\right)
.%
\end{array}%
\right.  \label{4.14}
\end{equation}%
Next, 
\begin{equation}
v\left( x,T/2\right) =v\left( x,t\right) -\dint\limits_{T/2}^{t}v_{t}\left(
x,\tau \right) d\tau .  \label{4.140}
\end{equation}%
Substituting this in (\ref{4.12}), we obtain%
\begin{equation}
\widetilde{k}\left( x\right) =2\left( \nabla u_{0,1}\left( x\right) \right)
^{-2}\left( v\left( x,t\right) -\dint\limits_{T/2}^{t}v_{t}\left( x,\tau
\right) d\tau \right) +F\left( x\right) .  \label{4.15}
\end{equation}%
Differentiating equations (\ref{4.9}) and (\ref{4.90}) with respect to $t$,
and using (\ref{4.10})-(\ref{4.15}), we obtain two integral differential
equations. The first equation is:%
\begin{equation}
\left. 
\begin{array}{c}
v_{t}+\Delta v+\dint\limits_{\Omega }Y\left( x,y\right) q\left( y,t\right)
dy+fq+f_{t}\dint\limits_{T/2}^{t}q\left( x,\tau \right) d\tau - \\ 
-k_{2}\nabla v\left( \nabla u_{1}+\nabla u_{2}\right) /2-k_{2}\left(
\dint\limits_{T/2}^{t}\nabla v\left( x,\tau \right) d\tau \right) \left(
\nabla u_{1t}+\nabla u_{2t}\right) /2 \\ 
-\left( \nabla u_{0,1}\left( x\right) \right) ^{-2}\left( \nabla
u_{1}\right) ^{2}\left( v\left( x,t\right)
-\dint\limits_{T/2}^{t}v_{t}\left( x,\tau \right) d\tau \right) = \\ 
=F\left( x\right) \left( \nabla u_{1}\right) ^{2}/2-f_{t}\widetilde{m}%
_{0}+k_{2}\text{ }\nabla \widetilde{u}_{0}\left( \nabla u_{1t}+\nabla
u_{2t}\right) /2,\text{ }\left( x,t\right) \in Q_{T}.%
\end{array}%
\right.  \label{4.16}
\end{equation}%
And the second equation is 
\begin{equation}
\left. 
\begin{array}{c}
q_{t}-\Delta q-\func{div}\left( k_{2}q\nabla u_{1}\right) -\func{div}\left(
k_{2}\nabla u_{1t}\dint\limits_{T/2}^{t}q\left( x,\tau \right) d\tau \right)
- \\ 
-\func{div}\left( k_{2}m_{2}\nabla v\right) -\func{div}\left(
k_{2}m_{2t}\dint\limits_{T/2}^{t}\nabla v\left( x,\tau \right) d\tau \right)
- \\ 
-\func{div}\left( \left( 2\left( \nabla u_{0,1}\right) ^{-2}\left( v\left(
x,t\right) -\dint\limits_{T/2}^{t}v_{t}\left( x,\tau \right) d\tau \right)
\right) m_{1}\nabla u_{1}\right) = \\ 
=\func{div}\left( Fm_{1}\nabla u_{1}\right) +\func{div}\left( k_{2}\nabla
u_{1t}\widetilde{m}_{0}\right) +\func{div}\left( k_{2}m_{2t}\nabla 
\widetilde{u}_{0}\right) ,\text{ }\left( x,t\right) \in Q_{T}.%
\end{array}%
\right.  \label{4.17}
\end{equation}%
Since the term 
\begin{equation}
\dint\limits_{T/2}^{t}v_{t}\left( x,\tau \right) d\tau  \label{4.170}
\end{equation}%
is present in both equations (\ref{4.16}) and (\ref{4.17}), then we
differentiate equations (\ref{4.16}) and (\ref{4.17}) with respect to $t$
once again. Denote 
\begin{equation}
w\left( x,t\right) =v_{t}\left( x,t\right) ,\text{ }r\left( x,t\right)
=q_{t}\left( x,t\right) .  \label{4.18}
\end{equation}%
Then term (\ref{4.170}) is replaced in the resulting equations either with $%
v\left( x,t\right) $ or with%
\begin{equation}
\dint\limits_{T/2}^{t}w\left( x,\tau \right) d\tau .  \label{4.180}
\end{equation}%
Using (\ref{4.16}), (\ref{4.18}) and (\ref{4.180}), we obtain 
\begin{equation}
\left. 
\begin{array}{c}
w_{t}+\Delta w+\dint\limits_{\Omega }Y\left( x,y\right) r\left( y,t\right)
dy+2f_{t}q+fr+f_{tt}\dint\limits_{T/2}^{t}q\left( x,\tau \right) d\tau - \\ 
-k_{2}\nabla w\left( \nabla u_{1}+\nabla u_{2}\right) /2-k_{2}\nabla v\left(
\nabla u_{1t}+\nabla u_{2t}\right) - \\ 
-k_{2}\left( \dint\limits_{T/2}^{t}\nabla v\left( x,\tau \right) d\tau
\right) \left( \nabla u_{1tt}+\nabla u_{2tt}\right) /2- \\ 
-2\left( \nabla u_{0,1}\left( x\right) \right) ^{-2}\nabla u_{1t}\nabla
u_{1}\left( v\left( x,t\right) -\dint\limits_{T/2}^{t}w\left( x,\tau \right)
d\tau \right) = \\ 
=F\left( x\right) \left( \nabla u_{1t}\nabla u_{1}\right) -f_{tt}\widetilde{m%
}_{0}+k_{2}\text{ }\nabla \widetilde{u}_{0}\left( \nabla u_{1tt}+\nabla
u_{2tt}\right) /2,\text{ }\left( x,t\right) \in Q_{T}.%
\end{array}%
\right.  \label{4.19}
\end{equation}%
Also, (\ref{4.17}) and (\ref{4.18}) lead to%
\begin{equation}
\left. 
\begin{array}{c}
r_{t}-\Delta r-2\func{div}\left( k_{2}r\nabla u_{1}\right) -\func{div}\left(
k_{2}q\nabla u_{1t}\right) - \\ 
-\func{div}\left( k_{2}m_{2}\nabla w\right) -2\func{div}\left(
k_{2}m_{2t}\nabla v\right) -\func{div}\left(
k_{2}m_{2tt}\dint\limits_{T/2}^{t}\nabla v\left( x,\tau \right) d\tau
\right) - \\ 
-\func{div}\left( \left( -2\left( \nabla u_{0,1}\left( x\right) \right)
^{-2}\left( v\left( x,t\right) -\dint\limits_{T/2}^{t}w\left( x,\tau \right)
d\tau \right) \right) \left( m_{1}\nabla u_{1}\right) _{t}\right) = \\ 
=\func{div}\left( F\left( x\right) \left( m_{1}\nabla u_{1}\right)
_{t}\right) +\func{div}\left( k_{2}\nabla u_{1tt}\widetilde{m}_{0}\right) +
\\ 
+\func{div}\left( k_{2}m_{2tt}\nabla \widetilde{u}_{0}\right) ,\text{ }%
\left( x,t\right) \in Q_{T}.%
\end{array}%
\right.  \label{4.20}
\end{equation}

It is well known that Carleman estimates can work not only with Partial
Differential Equations but with differential inequalities as well. Hence, it
is convenient to rewrite equations (\ref{4.16}), (\ref{4.17}), (\ref{4.19})
and (\ref{4.20}) in the more compact forms of differential inequalities.
Thus, using (\ref{2.2}), (\ref{4.1})-(\ref{4.7}), (\ref{4.14}) and (\ref%
{4.16}), we obtain: 
\begin{equation}
\left. 
\begin{array}{c}
\left\vert v_{t}+\Delta v\right\vert \leq C\left( \left\vert \nabla
v\right\vert +\left\vert v\right\vert +\left\vert
\dint\limits_{T/2}^{t}\left\vert \nabla v\right\vert \left( x,\tau \right)
d\tau \right\vert +\left\vert \dint\limits_{T/2}^{t}\left\vert w\right\vert
\left( x,\tau \right) d\tau \right\vert \right) + \\ 
+C\left( \left\vert \dint\limits_{T/2}^{t}\left\vert q\right\vert \left(
x,\tau \right) d\tau \right\vert +G\left( q\right) \left( x,t\right)
+\left\vert q\right\vert \right) +C\overline{\delta }\left( x,t\right) ,%
\text{ }\left( x,t\right) \in Q_{T}.%
\end{array}%
\right.  \label{4.21}
\end{equation}%
In (\ref{4.21})%
\begin{equation}
\left. 
\begin{array}{c}
G\left( q\right) \left( x,t\right) = \\ 
=\left\{ 
\begin{array}{c}
\dint\limits_{\Omega _{1}}\left\vert q\right\vert \left( x_{1},\overline{y}%
,t\right) d\overline{y}\text{ in the case } \\ 
\text{of the second line of (\ref{4.800}),} \\ 
\dint\limits_{\Omega _{1}}\left( \dint\limits_{x_{1}}^{b}\left\vert
q\right\vert \left( y_{1},\overline{y},t\right) dy_{1}\right) d\overline{y}%
\text{ in the case of } \\ 
\text{the third line of (\ref{4.800}).}%
\end{array}%
\right.%
\end{array}%
\right.  \label{4.210}
\end{equation}%
Using (\ref{3.7}), (\ref{3.8}) and (\ref{4.210}), we obtain%
\begin{equation}
\dint\limits_{Q_{T}}\left[ G\left( h\right) \left( x,t\right) \right]
^{2}\varphi _{\lambda }\left( x,t\right) dxdt\leq \widetilde{C}%
\dint\limits_{Q_{T}}h^{2}\varphi _{\lambda }dxdt,\text{ }\forall h\in
L_{2}\left( Q_{T}\right) ,\text{ }\forall \lambda >0.  \label{4.211}
\end{equation}%
Using (\ref{4.17}), we also obtain similarly with the above%
\begin{equation}
\left. 
\begin{array}{c}
\left\vert q_{t}-\Delta q\right\vert \leq C\left( \left\vert \nabla
q\right\vert +\left\vert q\right\vert +\left\vert
\dint\limits_{T/2}^{t}\left( \left\vert \nabla q\right\vert +\left\vert
q\right\vert \right) \left( x,\tau \right) d\tau \right\vert \right) + \\ 
+C\left( \left\vert \Delta v\right\vert +\left\vert \nabla v\right\vert
+\left\vert \dint\limits_{T/2}^{t}\left( \left\vert \Delta v\right\vert
+\left\vert \nabla v\right\vert \right) \left( x,\tau \right) d\tau
\right\vert \right) + \\ 
+C\left\vert \dint\limits_{T/2}^{t}\left( \left\vert \nabla w\right\vert
+\left\vert w\right\vert \right) \left( x,\tau \right) d\tau \right\vert +C%
\overline{\delta }\left( x,t\right) ,\text{ }\left( x,t\right) \in Q_{T}.%
\end{array}%
\right.  \label{4.22}
\end{equation}%
Similarly, (\ref{4.18}) and (\ref{4.19}) lead to:%
\begin{equation}
\left. 
\begin{array}{c}
\left\vert w_{t}+\Delta w\right\vert \leq C\left( \left\vert \nabla
w\right\vert +\left\vert w\right\vert +\left\vert
\dint\limits_{T/2}^{t}\left\vert w\right\vert \left( x,\tau \right) d\tau
\right\vert \right) + \\ 
+C\left( \left\vert \nabla v\right\vert +\left\vert v\right\vert +\left\vert
\dint\limits_{T/2}^{t}\left\vert \nabla v\right\vert \left( x,\tau \right)
d\tau \right\vert \right) +C\left\vert r\right\vert + \\ 
+C\left( \left\vert q\right\vert +\left\vert
\dint\limits_{T/2}^{t}\left\vert q\right\vert \left( x,\tau \right) d\tau
\right\vert +G\left( r\right) \left( x,t\right) \right) +C\overline{\delta }%
\left( x,t\right) ,\text{ }\left( x,t\right) \in Q_{T}.%
\end{array}%
\right.  \label{4.23}
\end{equation}%
Finally, using (\ref{4.1})-(\ref{4.7}), (\ref{4.14}) and (\ref{4.20}), we
obtain%
\begin{equation}
\left. 
\begin{array}{c}
\left\vert r_{t}-\Delta r\right\vert \leq C\left( \left\vert \nabla
r\right\vert +\left\vert r\right\vert +\left\vert \nabla q\right\vert
+\left\vert q\right\vert +\left\vert \Delta v\right\vert +\left\vert \nabla
v\right\vert +\left\vert v\right\vert \right) + \\ 
+C\left( \left\vert \Delta w\right\vert +\left\vert \nabla w\right\vert
+\left\vert w\right\vert \right) + \\ 
+C\left\vert \dint\limits_{T/2}^{t}\left( \left\vert \Delta v\right\vert
+\left\vert \nabla v\right\vert +\left\vert v\right\vert \right) \left(
x,\tau \right) d\tau \right\vert + \\ 
+C\left\vert \dint\limits_{T/2}^{t}\left( \left\vert \nabla w\right\vert
+\left\vert w\right\vert \right) \left( x,\tau \right) d\tau \right\vert +C%
\overline{\delta }\left( x,t\right) ,\text{ }\left( x,t\right) \in Q_{T}.%
\end{array}%
\right.  \label{4.24}
\end{equation}%
In (\ref{4.21})-(\ref{4.24}) $\overline{\delta }\left( x,t\right) \in
L_{2}\left( Q_{T}\right) $ denotes different non-negative functions such that%
\begin{equation}
\left\Vert \overline{\delta }\left( x,t\right) \right\Vert _{L_{2}\left(
Q_{T}\right) }\leq \delta .  \label{4.240}
\end{equation}%
Using (\ref{4.5}), (\ref{4.10}) and (\ref{4.18}) we obtain the following
lateral Cauchy data for the system of inequalities (\ref{4.21})-(\ref{4.24}):%
\begin{equation}
\left. 
\begin{array}{c}
v\mid _{S_{T}}=\widetilde{g}_{0t}\left( x,t\right) ,\text{ }\partial
_{n}v\mid _{S_{T}}=\widetilde{g}_{1t}\left( x,t\right) , \\ 
q\mid _{S_{T}}=\widetilde{p}_{0t}\left( x,t\right) ,\text{ }\partial
_{n}q\mid _{S_{T}}=\widetilde{p}_{1t}\left( x,t\right) , \\ 
w\mid _{S_{T}}=\widetilde{g}_{0tt}\left( x,t\right) ,\text{ }\partial
_{n}w\mid _{S_{T}}=\widetilde{g}_{1tt}\left( x,t\right) , \\ 
r\mid _{S_{T}}=\widetilde{p}_{0tt}\left( x,t\right) ,\text{ }\partial
_{n}r\mid _{S_{T}}=\widetilde{p}_{1tt}\left( x,t\right) .%
\end{array}%
\right.  \label{4.25}
\end{equation}

Square both sides of each of inequalities (\ref{4.21})-(\ref{4.24}). Then
multiply each of the resulting inequalities by the Carleman Weight Functions 
$\varphi _{\lambda }\left( x_{1},t\right) $ defined in (\ref{3.01}) and
apply the Carleman estimates (\ref{3.04}) using the first line of (\ref{3.2}%
), where $\lambda _{0}\geq 1$ was chosen in Theorem 3.1. Also, use second
and third lines of (\ref{4.7}) as well as (\ref{4.240}) and (\ref{4.25}). In
addition, use use (\ref{3.10}), (\ref{4.210}) and (\ref{4.211}). We obtain
four inequalities.

\subsection{Estimating functions $v$ and $q$}

\label{sec:4.1}

The first inequality is: 
\begin{equation}
\left. 
\begin{array}{c}
\left( 1/\lambda \right) \dint\limits_{Q_{T}}\left(
v_{t}^{2}+\dsum\limits_{i,j=1}^{n}v_{x_{i}x_{j}}^{2}\right) \varphi
_{\lambda }dxdt+\dint\limits_{Q_{T}}\left( \lambda \left( \nabla v\right)
^{2}+\lambda ^{3}v^{2}\right) \varphi _{\lambda }dxdt \\ 
\leq C\text{ }\dint\limits_{Q_{T}}\left[ \left\vert \nabla v\right\vert
^{2}+\left\vert v\right\vert ^{2}\right] \varphi _{\lambda }dxdt+ \\ 
+C\text{ }\dint\limits_{Q_{T}}w^{2}\varphi _{\lambda
}dxdt+C\dint\limits_{Q_{T}}q^{2}\varphi _{\lambda }dxdt+ \\ 
+C\left( \left\Vert v\left( x,T\right) \right\Vert _{H^{1}\left( \Omega
\right) }^{2}+\left\Vert v\left( x,0\right) \right\Vert _{H^{1}\left( \Omega
\right) }^{2}\right) \exp \left( -2\lambda \left( \alpha
T^{2}/4-b^{2}\right) \right) + \\ 
+Ce^{3\lambda b^{2}}\delta ^{2},\text{ }\forall \lambda \geq \lambda _{0}.%
\end{array}%
\right.  \label{4.26}
\end{equation}

Choose a sufficiently large number 
\begin{equation}
\lambda _{1}=\lambda _{1}\left( N,\varepsilon ,\Omega ,T,c\right) \geq
\lambda _{0}\geq 1  \label{4.280}
\end{equation}%
such that $\lambda _{1}>2C.$ Then (\ref{4.26}) and trace theorem lead to%
\begin{equation}
\left. 
\begin{array}{c}
\left( 1/\lambda \right) \dint\limits_{Q_{T}}\left(
v_{t}^{2}+\dsum\limits_{i,j=1}^{n}v_{x_{i}x_{j}}^{2}\right) \varphi
_{\lambda }dxdt+\dint\limits_{Q_{T}}\left( \lambda \left( \nabla v\right)
^{2}+\lambda ^{3}v^{2}\right) \varphi _{\lambda }dxdt\leq \\ 
\leq C\dint\limits_{Q_{T}}\left( w^{2}+q^{2}\right) \varphi _{\lambda }dxdt+
\\ 
+C\exp \left( -2\lambda \left( \alpha T^{2}/4-b^{2}\right) \right)
\left\Vert v\right\Vert _{H^{2}\left( Q_{T}\right) }^{2}+C\delta
^{2}e^{3\lambda b^{2}},\text{ }\forall \lambda \geq \lambda _{1}.%
\end{array}%
\right.  \label{4.29}
\end{equation}%
Applying similar arguments to (\ref{4.22}), we obtain 
\begin{equation}
\left. 
\begin{array}{c}
\left( 1/\lambda \right) \dint\limits_{Q_{T}}\left(
q_{t}^{2}+\dsum\limits_{i,j=1}^{n}q_{x_{i}x_{j}}^{2}\right) \varphi
_{\lambda }dxdt+\dint\limits_{Q_{T}}\left( \lambda \left( \nabla q\right)
^{2}+\lambda ^{3}q^{2}\right) \varphi _{\lambda }dxdt\leq \\ 
+C\dint\limits_{Q_{T}}\left( \left( \Delta v\right) ^{2}+\left( \nabla
v\right) ^{2}+\left( \nabla w\right) ^{2}+w^{2}\right) \varphi _{\lambda
}dxdt+ \\ 
+C\exp \left( -2\lambda \left( \alpha T^{2}/4-b^{2}\right) \right)
\left\Vert q\right\Vert _{H^{2}\left( Q_{T}\right) }^{2}+C\delta
^{2}e^{3\lambda b^{2}},\text{ }\forall \lambda \geq \lambda _{1}.%
\end{array}%
\right.  \label{4.32}
\end{equation}%
It follows from (\ref{4.32}) that%
\begin{equation}
\left. 
\begin{array}{c}
\dint\limits_{Q_{T}}q^{2}\varphi _{\lambda }dxdt\leq \left( C/\lambda
^{3}\right) \dint\limits_{Q_{T}}\left( \left( \Delta v\right) ^{2}+\left(
\nabla v\right) ^{2}+\left( \nabla w\right) ^{2}+w^{2}\right) \varphi
_{\lambda }dxdt+ \\ 
+C\exp \left( -2\lambda \left( \alpha T^{2}/4-b^{2}\right) \right)
\left\Vert q\right\Vert _{H^{2}\left( Q_{T}\right) }^{2}+C\delta
^{2}e^{3\lambda b^{2}},\text{ }\forall \lambda \geq \lambda _{1}.%
\end{array}%
\right.  \label{4.33}
\end{equation}%
Substituting (\ref{4.33}) in (\ref{4.29}) and using $\left( 1/\lambda
\right) >>\left( 1/\lambda ^{3}\right) $ for all $\lambda \geq \lambda _{1},$
we obtain 
\begin{equation}
\left. 
\begin{array}{c}
\left( 1/\lambda \right) \dint\limits_{Q_{T}}\left(
v_{t}^{2}+\dsum\limits_{i,j=1}^{n}v_{x_{i}x_{j}}^{2}\right) \varphi
_{\lambda }dxdt+\dint\limits_{Q_{T}}\left( \lambda \left( \nabla v\right)
^{2}+\lambda ^{3}v^{2}\right) \varphi _{\lambda }dxdt\leq \\ 
\leq C\dint\limits_{Q_{T}}\left[ \left( 1/\lambda ^{3}\right) \left( \nabla
w\right) ^{2}+w^{2}\right] \varphi _{\lambda }dxdt+ \\ 
+C\exp \left( -2\lambda \left( \alpha T^{2}/4-b^{2}\right) \right) \left(
\left\Vert v\right\Vert _{H^{2}\left( Q_{T}\right) }^{2}+\left\Vert
q\right\Vert _{H^{2}\left( Q_{T}\right) }^{2}\right) + \\ 
+C\delta ^{2}e^{3\lambda b^{2}},\text{ }\forall \lambda \geq \lambda _{1}.%
\end{array}%
\right.  \label{4.34}
\end{equation}%
In particular, (\ref{4.34}) implies%
\begin{equation}
\left. 
\begin{array}{c}
\dint\limits_{Q_{T}}\left( \Delta v\right) ^{2}\varphi _{\lambda }dxdt\leq
C\dint\limits_{Q_{T}}\left[ \left( 1/\lambda ^{2}\right) \left( \nabla
w\right) ^{2}+\lambda w^{2}\right] \varphi _{\lambda }dxdt+ \\ 
+C\lambda \exp \left( -2\lambda \left( \alpha T^{2}/4-b^{2}\right) \right)
\left( \left\Vert v\right\Vert _{H^{2}\left( Q_{T}\right) }^{2}+\left\Vert
q\right\Vert _{H^{2}\left( Q_{T}\right) }^{2}\right) + \\ 
+C\delta ^{2}\lambda e^{3\lambda b^{2}},\text{ }\forall \lambda \geq \lambda
_{1}.%
\end{array}%
\right.  \label{4.35}
\end{equation}%
Substituting (\ref{4.35}) in (\ref{4.32}), we obtain%
\begin{equation}
\left. 
\begin{array}{c}
\left( 1/\lambda \right) \dint\limits_{Q_{T}}\left(
q_{t}^{2}+\dsum\limits_{i,j=1}^{n}q_{x_{i}x_{j}}^{2}\right) \varphi
_{\lambda }dxdt+\dint\limits_{Q_{T}}\left( \lambda \left( \nabla q\right)
^{2}+\lambda ^{3}q^{2}\right) \varphi _{\lambda }dxdt\leq \\ 
+C\dint\limits_{Q_{T}}\left( \left( \nabla v\right) ^{2}+\left( \nabla
w\right) ^{2}+\lambda w^{2}\right) \varphi _{\lambda }dxdt+ \\ 
+C\exp \left( -2\lambda \left( \alpha T^{2}/4-b^{2}\right) \right)
\left\Vert q\right\Vert _{H^{2}\left( Q_{T}\right) }^{2}+C\delta
^{2}e^{3\lambda b^{2}},\text{ }\forall \lambda \geq \lambda _{1}.%
\end{array}%
\right.  \label{4.350}
\end{equation}%
Summing up (\ref{4.34}) and (\ref{4.350}), we obtain%
\begin{equation}
\left. 
\begin{array}{c}
\left( 1/\lambda \right) \dint\limits_{Q_{T}}\left[ \left(
v_{t}^{2}+q_{t}^{2}\right) +\dsum\limits_{i,j=1}^{n}\left(
v_{x_{i}x_{j}}^{2}+q_{x_{i}x_{j}}^{2}\right) \right] \varphi _{\lambda }dxdt+
\\ 
+\dint\limits_{Q_{T}}\left[ \lambda \left( \left( \nabla v\right)
^{2}+\left( \nabla q\right) ^{2}\right) +\lambda ^{3}\left(
v^{2}+q^{2}\right) \right] \varphi _{\lambda }dxdt\leq \\ 
\leq C\dint\limits_{Q_{T}}\left( \left( \nabla w\right) ^{2}+\lambda
w^{2}\right) \varphi _{\lambda }dxdt+ \\ 
+C\lambda \exp \left( -2\lambda \left( \alpha T^{2}/4-b^{2}\right) \right)
\left( \left\Vert v\right\Vert _{H^{2}\left( Q_{T}\right) }^{2}+\left\Vert
q\right\Vert _{H^{2}\left( Q_{T}\right) }^{2}\right) + \\ 
+C\delta ^{2}\lambda e^{3\lambda b^{2}},\text{ }\forall \lambda \geq \lambda
_{1}.%
\end{array}%
\right.  \label{4.36}
\end{equation}

\subsection{\ Estimating functions $w$ and $r$}

\label{sec:4.2}

We \ now estimate functions $w$ and $r$. Using (\ref{4.23}), we obtain
similarly with (\ref{4.29})%
\begin{equation}
\left. 
\begin{array}{c}
\left( 1/\lambda \right) \dint\limits_{Q_{T}}\left(
w_{t}^{2}+\dsum\limits_{i,j=1}^{n}w_{x_{i}x_{j}}^{2}\right) \varphi
_{\lambda }dxdt+\dint\limits_{Q_{T}}\left( \lambda \left( \nabla w\right)
^{2}+\lambda ^{3}w^{2}\right) \varphi _{\lambda }dxdt\leq \\ 
\leq C\dint\limits_{Q_{T}}\left( \left( \nabla v\right)
^{2}+v^{2}+q^{2}\right) \varphi _{\lambda
}dxdt+C\dint\limits_{Q_{T}}r^{2}\varphi _{\lambda }dxdt \\ 
+C\exp \left( -2\lambda \left( \alpha T^{2}/4-b^{2}\right) \right)
\left\Vert w\right\Vert _{H^{2}\left( Q_{T}\right) }^{2}+C\delta
^{2}e^{3\lambda b^{2}},\text{ }\forall \lambda \geq \lambda _{1}.%
\end{array}%
\right.  \label{4.37}
\end{equation}%
Next, using (\ref{4.24}), we obtain similarly with (\ref{4.32})%
\begin{equation}
\left. 
\begin{array}{c}
\left( 1/\lambda \right) \dint\limits_{Q_{T}}\left(
r_{t}^{2}+\dsum\limits_{i,j=1}^{n}r_{x_{i}x_{j}}^{2}\right) \varphi
_{\lambda }dxdt+\dint\limits_{Q_{T}}\left( \lambda \left( \nabla r\right)
^{2}+\lambda ^{3}r^{2}\right) \varphi _{\lambda }dxdt\leq \\ 
+C\dint\limits_{Q_{T}}\left( \left( \Delta v\right) ^{2}+\left( \nabla
v\right) ^{2}+v^{2}\right) \varphi _{\lambda }dxdt+ \\ 
+C\dint\limits_{Q_{T}}\left( \left( \Delta w\right) ^{2}+\left( \nabla
w\right) ^{2}+w^{2}\right) \varphi _{\lambda }dxdt+ \\ 
+C\exp \left( -2\lambda \left( \alpha T^{2}/4-b^{2}\right) \right)
\left\Vert r\right\Vert _{H^{2}\left( Q_{T}\right) }^{2}+C\delta
^{2}e^{3\lambda b^{2}},\text{ }\forall \lambda \geq \lambda _{1}.%
\end{array}%
\right.  \label{4.38}
\end{equation}%
By (\ref{4.35}) and (\ref{4.36})%
\begin{equation}
\left. 
\begin{array}{c}
\dint\limits_{Q_{T}}\left( \left( \Delta v\right) ^{2}+\left( \nabla
v\right) ^{2}+v^{2}\right) \varphi _{\lambda }dxdt\leq C\dint\limits_{Q_{T}} 
\left[ \left( 1/\lambda ^{2}\right) \left( \nabla w\right) ^{2}+\lambda w^{2}%
\right] \varphi _{\lambda }dxdt+ \\ 
+C\lambda \exp \left( -2\lambda \left( \alpha T^{2}/4-b^{2}\right) \right)
\left( \left\Vert v\right\Vert _{H^{2}\left( Q_{T}\right) }^{2}+\left\Vert
q\right\Vert _{H^{2}\left( Q_{T}\right) }^{2}\right) + \\ 
+C\delta ^{2}\lambda e^{3\lambda b^{2}},\text{ }\forall \lambda \geq \lambda
_{1}.%
\end{array}%
\right.  \label{4.39}
\end{equation}%
Substituting (\ref{4.39}) in (\ref{4.38}), we obtain%
\begin{equation}
\left. 
\begin{array}{c}
\left( 1/\lambda \right) \dint\limits_{Q_{T}}\left(
r_{t}^{2}+\dsum\limits_{i,j=1}^{n}r_{x_{i}x_{j}}^{2}\right) \varphi
_{\lambda }dxdt+\dint\limits_{Q_{T}}\left( \lambda \left( \nabla r\right)
^{2}+\lambda ^{3}r^{2}\right) \varphi _{\lambda }dxdt\leq \\ 
\leq C\dint\limits_{Q_{T}}\left( \left( \Delta w\right) ^{2}+\left( \nabla
w\right) ^{2}+\lambda w^{2}\right) \varphi _{\lambda }dxdt+ \\ 
+C\lambda \exp \left( -2\lambda \left( \alpha T^{2}/4-b^{2}\right) \right)
\left( \left\Vert v\right\Vert _{H^{2}\left( Q_{T}\right) }^{2}+\left\Vert
q\right\Vert _{H^{2}\left( Q_{T}\right) }^{2}+\left\Vert r\right\Vert
_{H^{2}\left( Q_{T}\right) }^{2}\right) \\ 
+C\delta ^{2}\lambda e^{3\lambda b^{2}},\text{ }\forall \lambda \geq \lambda
_{1}.%
\end{array}%
\right.  \label{4.40}
\end{equation}%
Next, substituting (\ref{4.39}) in (\ref{4.37}), we obtain 
\begin{equation}
\left. 
\begin{array}{c}
\left( 1/\lambda \right) \dint\limits_{Q_{T}}\left(
w_{t}^{2}+\dsum\limits_{i,j=1}^{n}w_{x_{i}x_{j}}^{2}\right) \varphi
_{\lambda }dxdt+\dint\limits_{Q_{T}}\left( \lambda \left( \nabla w\right)
^{2}+\lambda ^{3}w^{2}\right) \varphi _{\lambda }dxdt\leq \\ 
\leq C\dint\limits_{Q_{T}}\left( q^{2}+r^{2}\right) \varphi _{\lambda }dxdt+
\\ 
+C\lambda \exp \left( -2\lambda \left( \alpha T^{2}/4-b^{2}\right) \right)
\times \\ 
\times \left( \left\Vert v\right\Vert _{H^{2}\left( Q_{T}\right)
}^{2}+\left\Vert q\right\Vert _{H^{2}\left( Q_{T}\right) }^{2}+\left\Vert
w\right\Vert _{H^{2}\left( Q_{T}\right) }^{2}+\left\Vert r\right\Vert
_{H^{2}\left( Q_{T}\right) }^{2}\right) + \\ 
+C\delta ^{2}\lambda e^{3\lambda b^{2}},\text{ }\forall \lambda \geq \lambda
_{1}.%
\end{array}%
\right.  \label{4.41}
\end{equation}%
In particular, (\ref{4.40}) implies%
\begin{equation}
\left. 
\begin{array}{c}
\dint\limits_{Q_{T}}r^{2}\varphi _{\lambda }dxdt\leq \left( C/\lambda
^{3}\right) \dint\limits_{Q_{T}}\left( \left( \Delta w\right) ^{2}+\left(
\nabla w\right) ^{2}\right) \varphi _{\lambda }dxdt+ \\ 
+\left( C/\lambda ^{2}\right) \dint\limits_{Q_{T}}w^{2}\varphi _{\lambda
}dxdt+] \\ 
+C\exp \left( -2\lambda \left( \alpha T^{2}/4-b^{2}\right) \right) \left(
\left\Vert v\right\Vert _{H^{2}\left( Q_{T}\right) }^{2}+\left\Vert
q\right\Vert _{H^{2}\left( Q_{T}\right) }^{2}+\left\Vert r\right\Vert
_{H^{2}\left( Q_{T}\right) }^{2}\right) + \\ 
+C\delta ^{2}e^{3\lambda b^{2}},\text{ }\forall \lambda \geq \lambda _{1}.%
\end{array}%
\right.  \label{4.42}
\end{equation}%
Substituting (\ref{4.42}) in (\ref{4.41}), we obtain%
\begin{equation}
\left. 
\begin{array}{c}
\left( 1/\lambda \right) \dint\limits_{Q_{T}}\left(
w_{t}^{2}+\dsum\limits_{i,j=1}^{n}w_{x_{i}x_{j}}^{2}\right) \varphi
_{\lambda }dxdt+\dint\limits_{Q_{T}}\left( \lambda \left( \nabla w\right)
^{2}+\lambda ^{3}w^{2}\right) \varphi _{\lambda }dxdt\leq \\ 
\leq C\dint\limits_{Q_{T}}q^{2}\varphi _{\lambda }dxdt+ \\ 
+C\lambda \exp \left( -2\lambda \left( \alpha T^{2}/4-b^{2}\right) \right)
\times \\ 
\times \left( \left\Vert v\right\Vert _{H^{2}\left( Q_{T}\right)
}^{2}+\left\Vert q\right\Vert _{H^{2}\left( Q_{T}\right) }^{2}+\left\Vert
w\right\Vert _{H^{2}\left( Q_{T}\right) }^{2}+\left\Vert r\right\Vert
_{H^{2}\left( Q_{T}\right) }^{2}\right) + \\ 
+C\delta ^{2}\lambda e^{3\lambda b^{2}},\text{ }\forall \lambda \geq \lambda
_{1}.%
\end{array}%
\right.  \label{4.43}
\end{equation}%
In particular, (\ref{4.43}) implies%
\begin{equation}
\left. 
\begin{array}{c}
\dint\limits_{Q_{T}}\left( \Delta w\right) ^{2}\varphi _{\lambda }dxdt\leq
C\lambda \dint\limits_{Q_{T}}q^{2}\varphi _{\lambda }dxdt+ \\ 
+C\lambda ^{2}\exp \left( -2\lambda \left( \alpha T^{2}/4-b^{2}\right)
\right) \times \\ 
\times \left( \left\Vert v\right\Vert _{H^{2}\left( Q_{T}\right)
}^{2}+\left\Vert q\right\Vert _{H^{2}\left( Q_{T}\right) }^{2}+\left\Vert
w\right\Vert _{H^{2}\left( Q_{T}\right) }^{2}+\left\Vert r\right\Vert
_{H^{2}\left( Q_{T}\right) }^{2}\right) + \\ 
+C\delta ^{2}\lambda ^{2}e^{3\lambda b^{2}},\text{ }\forall \lambda \geq
\lambda _{1}.%
\end{array}%
\right.  \label{4.44}
\end{equation}%
Substituting (\ref{4.44}) in (\ref{4.40}), we obtain%
\begin{equation}
\left. 
\begin{array}{c}
\left( 1/\lambda \right) \dint\limits_{Q_{T}}\left(
r_{t}^{2}+\dsum\limits_{i,j=1}^{n}r_{x_{i}x_{j}}^{2}\right) \varphi
_{\lambda }dxdt+\dint\limits_{Q_{T}}\left( \lambda \left( \nabla r\right)
^{2}+\lambda ^{3}r^{2}\right) \varphi _{\lambda }dxdt\leq \\ 
\leq C\dint\limits_{Q_{T}}\left( \left( \nabla w\right) ^{2}+\lambda
w^{2}+\lambda q^{2}\right) \varphi _{\lambda }dxdt+ \\ 
+C\lambda ^{2}\exp \left( -2\lambda \left( \alpha T^{2}/4-b^{2}\right)
\right) \times \\ 
\times \left( \left\Vert v\right\Vert _{H^{2}\left( Q_{T}\right)
}^{2}+\left\Vert q\right\Vert _{H^{2}\left( Q_{T}\right) }^{2}+\left\Vert
w\right\Vert _{H^{2}\left( Q_{T}\right) }^{2}+\left\Vert r\right\Vert
_{H^{2}\left( Q_{T}\right) }^{2}\right) + \\ 
+C\delta ^{2}\lambda ^{2}e^{3\lambda b^{2}},\text{ }\forall \lambda \geq
\lambda _{1}.%
\end{array}%
\right.  \label{4.45}
\end{equation}

Summing up (\ref{4.43}) and (\ref{4.45}), we obtain the following analog of (%
\ref{4.36}):%
\begin{equation}
\left. 
\begin{array}{c}
\left( 1/\lambda \right) \dint\limits_{Q_{T}}\left[ \left(
w_{t}^{2}+r_{t}^{2}\right) +\dsum\limits_{i,j=1}^{n}\left(
w_{x_{i}x_{j}}^{2}+r_{x_{i}x_{j}}^{2}\right) \right] \varphi _{\lambda }dxdt+
\\ 
+\dint\limits_{Q_{T}}\left[ \lambda \left( \left( \nabla w\right)
^{2}+\left( \nabla r\right) ^{2}\right) +\lambda ^{3}\left(
w^{2}+r^{2}\right) \right] \varphi _{\lambda }dxdt\leq C\lambda
\dint\limits_{Q_{T}}q^{2}\varphi _{\lambda }dxdt+ \\ 
+C\lambda ^{2}\exp \left( -2\lambda \left( \alpha T^{2}/4-b^{2}\right)
\right) \times \\ 
\times \left( \left\Vert v\right\Vert _{H^{2}\left( Q_{T}\right)
}^{2}+\left\Vert q\right\Vert _{H^{2}\left( Q_{T}\right) }^{2}+\left\Vert
w\right\Vert _{H^{2}\left( Q_{T}\right) }^{2}+\left\Vert r\right\Vert
_{H^{2}\left( Q_{T}\right) }^{2}\right) + \\ 
+C\delta ^{2}\lambda ^{2}e^{3\lambda b^{2}},\text{ }\forall \lambda \geq
\lambda _{1}.%
\end{array}%
\right.  \label{4.46}
\end{equation}

\subsection{Estimate all four functions $v,q,w,r$ simultaneously}

\label{sec:4.3}

Summing up (\ref{4.36}) and (\ref{4.46}), we obtain%
\begin{equation}
\left. 
\begin{array}{c}
\left( 1/\lambda \right) \dint\limits_{Q_{T}}\left[ \left(
v_{t}^{2}+q_{t}^{2}+w_{t}^{2}+r_{t}^{2}\right)
+\dsum\limits_{i,j=1}^{n}\left(
v_{x_{i}x_{j}}^{2}+q_{x_{i}x_{j}}^{2}+w_{x_{i}x_{j}}^{2}+r_{x_{i}x_{j}}^{2}%
\right) \right] \varphi _{\lambda }dxdt+ \\ 
+\lambda \dint\limits_{Q_{T}}\left[ \left( \left( \nabla v\right)
^{2}+\left( \nabla q\right) ^{2}+\left( \nabla w\right) ^{2}+\left( \nabla
r\right) ^{2}\right) \right] \varphi _{\lambda }dxdt+ \\ 
+\lambda ^{3}\dint\limits_{Q_{T}}\left( v^{2}+q^{2}+w^{2}+r^{2}\right)
\varphi _{\lambda }dxdt\leq \\ 
\leq C\lambda ^{2}\exp \left( -2\lambda \left( \alpha T^{2}/4-b^{2}\right)
\right) \times \\ 
\times \left( \left\Vert v\right\Vert _{H^{2}\left( Q_{T}\right)
}^{2}+\left\Vert q\right\Vert _{H^{2}\left( Q_{T}\right) }^{2}+\left\Vert
w\right\Vert _{H^{2}\left( Q_{T}\right) }^{2}+\left\Vert r\right\Vert
_{H^{2}\left( Q_{T}\right) }^{2}\right) + \\ 
+C\delta ^{2}\lambda ^{2}e^{3\lambda b^{2}},\text{ }\forall \lambda \geq
\lambda _{1}.%
\end{array}%
\right.  \label{4.47}
\end{equation}

\subsection{The target H\"{o}lder stability estimate}

\label{sec:4.4}

By (\ref{3.001}) and (\ref{3.2})%
\begin{equation}
\dint\limits_{Q_{T}}h^{2}\varphi _{\lambda }dxdt\geq \exp \left[ 2\lambda
\left( a^{2}-\alpha \left( T/2-\varepsilon \right) ^{2}\right) \right]
\dint\limits_{Q_{\varepsilon ,T}}h^{2}dxdt,\text{ }\forall h\in L_{2}\left(
Q_{T}\right) .  \label{4.48}
\end{equation}%
Introduce the vector function $V\left( x,t\right) $ as well as its two
norms, 
\begin{equation}
\left. 
\begin{array}{c}
V\left( x,t\right) =\left( v,q,,w,r\right) ^{T}\left( x,t\right) , \\ 
\left\Vert V\right\Vert _{H^{2,1}\left( Q_{\varepsilon ,T}\right)
}^{2}=\left\Vert v\right\Vert _{H^{2,1}\left( Q_{\varepsilon ,T}\right)
}^{2}+\left\Vert q\right\Vert _{H^{2,1}\left( Q_{\varepsilon ,T}\right)
}^{2}+\left\Vert w\right\Vert _{H^{2,1}\left( Q_{\varepsilon ,T}\right)
}^{2}+\left\Vert r\right\Vert _{H^{2,1}\left( Q_{\varepsilon ,T}\right)
}^{2}, \\ 
\left\Vert V\right\Vert _{H^{2}\left( Q_{T}\right) }^{2}=\left\Vert
v\right\Vert _{H^{2}\left( Q_{T}\right) }^{2}+\left\Vert q\right\Vert
_{H^{2}\left( Q_{T}\right) }^{2}+\left\Vert w\right\Vert _{H^{2}\left(
Q_{T}\right) }^{2}+\left\Vert r\right\Vert _{H^{2}\left( Q_{T}\right) }^{2}.%
\end{array}%
\right.  \label{4.49}
\end{equation}%
Obviously%
\begin{equation}
\left\Vert V\right\Vert _{H^{2}\left( Q_{T}\right) }\leq C.  \label{4.490}
\end{equation}%
Using (\ref{4.47})-(\ref{4.49}), we obtain%
\begin{equation}
\left. 
\begin{array}{c}
\exp \left[ 2\lambda \left( a^{2}-\alpha \left( T/2-\varepsilon \right)
^{2}\right) \right] \left\Vert V\right\Vert _{H^{2,1}\left( Q_{\varepsilon
,T}\right) }^{2}\leq \\ 
\leq C\lambda ^{2}\exp \left( -2\lambda \left( \alpha T^{2}/4-b^{2}\right)
\right) \left\Vert V\right\Vert _{H^{2}\left( Q_{T}\right) }^{2}+C\delta
^{2}\lambda ^{2}e^{3\lambda b^{2}},\text{ }\forall \lambda \geq \lambda _{1}.%
\end{array}%
\right.  \label{4.50}
\end{equation}%
Dividing both sides of (\ref{4.50}) by $\exp \left[ 2\lambda \left(
a^{2}-\alpha \left( T/2-\varepsilon \right) ^{2}\right) \right] $ and
noticing that $\lambda ^{2}e^{-2\lambda a^{2}}<1$ for $\lambda \geq \lambda
_{1},$ we obtain%
\begin{equation}
\left. 
\begin{array}{c}
\left\Vert V\right\Vert _{H^{2,1}\left( Q_{\varepsilon ,T}\right) }^{2}\leq
C\exp \left[ -2\lambda \left( \alpha \varepsilon \left( T-\varepsilon
\right) -b^{2}\right) \right] \left\Vert V\right\Vert _{H^{2}\left(
Q_{T}\right) }^{2}+ \\ 
+C\delta ^{2}\exp \left[ 2\lambda \left( 3b^{2}/2+\alpha \left(
T/2-\varepsilon \right) ^{2}\right) \right] ,\text{ }\forall \lambda \geq
\lambda _{1}.%
\end{array}%
\right.  \label{4.51}
\end{equation}%
Let the number $\beta >0.$ Choose the parameter $\alpha $ as%
\begin{equation}
\alpha =\frac{\left( 1+\beta \right) b^{2}}{\varepsilon \left( T-\varepsilon
\right) }.  \label{4.52}
\end{equation}%
Then in the first line of (\ref{4.51}) 
\begin{equation}
\exp \left[ -2\lambda \left( \alpha \varepsilon \left( T-\varepsilon \right)
-b^{2}\right) \right] =\exp \left( -2\lambda \beta b^{2}\right) .
\label{4.53}
\end{equation}%
Next, using (\ref{4.52}), we obtain in the second line of (\ref{4.51})%
\begin{equation}
\exp \left[ 2\lambda \left( 3b^{2}/2+\alpha \left( T/2-\varepsilon \right)
^{2}\right) \right] =\exp \left( 2\lambda d\right) ,  \label{4.54}
\end{equation}%
\begin{equation}
d=\left[ \frac{3}{2}+\frac{\left( 1+\beta \right) \left( T/2-\varepsilon
\right) ^{2}}{\left( \varepsilon \left( T-\varepsilon \right) \right) }%
\right] b^{2}.  \label{4.540}
\end{equation}

Recall that the number $\rho \in \left( 0,1\right) .$ Choose $\lambda
=\lambda \left( \delta \right) $ such that 
\begin{equation}
\delta ^{2}e^{2\lambda \left( \delta \right) d}=\delta ^{2-2\rho }.
\label{4.55}
\end{equation}%
Hence, in (\ref{4.54}) 
\begin{equation}
2\lambda \left( \delta \right) =\ln \left( \delta ^{-2\rho /d}\right) .
\label{4.56}
\end{equation}%
Since we must have 
\begin{equation}
\lambda \left( \delta \right) \geq \lambda _{1}=\lambda _{1}\left(
N,\varepsilon ,\Omega ,T,c\right) \geq 1,  \label{4.560}
\end{equation}%
see (\ref{4.280}), then we must have%
\begin{equation}
\delta \in \left( 0,\delta _{0}\right) ,\text{ }\delta _{0}=\delta
_{0}\left( N,\varepsilon ,\Omega ,T,c,\rho \right) =\exp \left( -\frac{%
\lambda _{1}d}{\rho }\right) \in \left( 0,1\right) .  \label{4.57}
\end{equation}%
Substituting $2\lambda \left( \delta \right) $ from (\ref{4.56}) in the
right hand side of (\ref{4.53}), we obtain 
\begin{equation}
\exp \left( -2\lambda \left( \delta \right) \beta b^{2}\right) =\delta
^{\left( 2\rho \beta b^{2}\right) /d}.  \label{4.58}
\end{equation}%
Next, using (\ref{4.540}), we obtain%
\begin{equation}
\frac{2\rho \beta b^{2}}{d}=\frac{2\rho \beta }{3/2+\left( 1+\beta \right)
\left( T/2-\varepsilon \right) ^{2}/\left( \varepsilon \left( T-\varepsilon
\right) \right) }.  \label{4.59}
\end{equation}%
By (\ref{4.54}) and (\ref{4.55}) the term in the second line of (\ref{4.51})
can be estimated from the above by $C\delta ^{2-2\rho }$ when $\delta
\rightarrow 0.$ Hence, we want the multiplier $C\exp \left[ -2\lambda \left(
\alpha \varepsilon \left( T-\varepsilon \right) -b^{2}\right) \right] $ in
the first line of (\ref{4.51}) not to exceed $C\delta ^{2-2\rho }$ as $%
\delta \rightarrow 0.$ Hence, it follows from (\ref{4.53}), (\ref{4.55}), (%
\ref{4.58}) and (\ref{4.59}) that we should have 
\begin{equation}
\frac{\rho \beta b^{2}}{d}=\frac{\rho \beta }{3/2+\left( 1+\beta \right)
\left( T/2-\varepsilon \right) ^{2}/\left( \varepsilon \left( T-\varepsilon
\right) \right) }\geq 1-\rho .  \label{4.60}
\end{equation}%
This is equivalent with%
\begin{equation}
\beta \left[ \rho -\left( 1-\rho \right) \frac{\left( T/2-\varepsilon
\right) ^{2}}{\varepsilon \left( T-\varepsilon \right) }\right] \geq \left(
1-\rho \right) \left[ \frac{3}{2}+\frac{\left( T/2-\varepsilon \right) ^{2}}{%
\varepsilon \left( T-\varepsilon \right) }\right] .  \label{4.61}
\end{equation}%
To find numbers $\beta >0$ satisfying (\ref{4.61}), it is necessary and
sufficient to have%
\begin{equation}
\rho -\left( 1-\rho \right) \frac{\left( T/2-\varepsilon \right) ^{2}}{%
\varepsilon \left( T-\varepsilon \right) }>0.  \label{4.62}
\end{equation}%
It obviously follows from (\ref{4.69}) that inequality (\ref{4.62}) holds.
Thus, if we choose any $\beta $ satisfying (\ref{4.61}), then (\ref{4.57})-(%
\ref{4.60}) imply 
\begin{equation*}
\exp \left( -2\lambda \left( \delta \right) \beta b^{2}\right) =\delta
^{\left( 2\rho \beta b^{2}\right) /d}\leq \delta ^{2-2\rho },\text{ }\forall
\delta \in \left( 0,\delta _{0}\right) .
\end{equation*}%
Combining this with (\ref{4.51})-(\ref{4.58}) and recalling (\ref{4.490}),
we obtain 
\begin{equation*}
\left\Vert V\right\Vert _{H^{2,1}\left( Q_{\varepsilon ,T}\right) }\leq
C\delta ^{1-\rho },\text{ }\forall \delta \in \left( 0,\delta _{0}\right) .
\end{equation*}%
Hence, by (\ref{4.49})%
\begin{equation}
\left\Vert v\right\Vert _{H^{2,1}\left( Q_{\varepsilon ,T}\right)
}+\left\Vert q\right\Vert _{H^{2,1}\left( Q_{\varepsilon ,T}\right)
}+\left\Vert w\right\Vert _{H^{2,1}\left( Q_{\varepsilon ,T}\right)
}+\left\Vert r\right\Vert _{H^{2,1}\left( Q_{\varepsilon ,T}\right) }\leq
C\delta ^{1-\rho },\text{ }\forall \delta \in \left( 0,\delta _{0}\right) .
\label{4.63}
\end{equation}%
Next, it follows from the first line of (\ref{4.80}), (\ref{4.10}), (\ref%
{4.11}), (\ref{4.18}) and (\ref{4.63}) that we also have:%
\begin{equation}
\left\Vert \partial _{t}^{s}\widetilde{u}\right\Vert _{H^{2,1}\left(
Q_{\varepsilon ,T}\right) },\left\Vert \partial _{t}^{s}\widetilde{m}%
\right\Vert _{H^{2,1}\left( Q_{\varepsilon ,T}\right) }\leq C\delta ^{1-\rho
},\text{ }\forall \delta \in \left( 0,\delta _{0}\right) ,\text{ }s=0,1,2.
\label{4.64}
\end{equation}

Finally, using the first line of (\ref{4.7}), (\ref{4.12})-(\ref{4.140}),
Lemma 3.3 and (\ref{4.63}), we obtain%
\begin{equation}
\left\Vert \widetilde{k}\right\Vert _{L_{2}\left( \Omega \right) }\leq
C\delta ^{1-\rho },\text{ }\forall \delta \in \left( 0,\delta _{0}\right) .
\label{4.65}
\end{equation}%
The target estimate (\ref{4.8}) of this theorem follows immediately from
estimates (\ref{4.64}) and (\ref{4.65}).

To obtain the uniqueness result claimed by this theorem, it is sufficient to
set in (\ref{4.7}) and (\ref{4.8}) $\delta =0.$ Then we obtain $\widetilde{u}%
\left( x,t\right) =\widetilde{m}\left( x,t\right) =0$ in $Q_{\varepsilon ,T}$
\ and%
\begin{equation}
\widetilde{k}\left( x\right) \equiv 0\text{ in }\Omega .  \label{4.66}
\end{equation}%
Next, applying Theorems 4.1 and 5.1 of \cite{MFG4} to the system (\ref{4.9}%
), (\ref{4.90}) with condition (\ref{4.66}) and zero Neumann boundary
condition at $S_{T}$ (see (\ref{4.5}), (\ref{4.7})) in domains $\Omega
\times \left( 0,\varepsilon \right) $ and $\Omega \times \left(
T-\varepsilon ,T\right) $ respectively, we obtain $\widetilde{u}\left(
x,t\right) \equiv \widetilde{m}\left( x,t\right) \equiv 0$ in $Q_{T}.$ \ \ $%
\square $

\section{Proof of Theorem 3.3}

\label{sec:5}

Estimate (\ref{4.26}) was obtained under the assumption that inequalities in
the second line of (\ref{4.7}). The same is true for the above estimates,
similar with (\ref{4.26}), for functions $q,w,r$. Now, however, it follows
from (\ref{4.80}) that we can use the implication of (\ref{3.06}) from (\ref%
{3.05}) and those estimates will be still valid. Then we can precisely
repeat the proof of Theorem 3.2 and obtain H\"{o}lder obtain estimates (\ref%
{4.8})\emph{.} $\ $

As to uniqueness, if $\delta =0,$ then the above paragraph implies of course
that $\widetilde{u}\left( x,t\right) =\widetilde{m}\left( x,t\right) =0$ in $%
Q_{\varepsilon ,T}$, and (\ref{4.66}) is also in place. As to Theorems 4.1
and 5.1 of \cite{MFG4} being applied to the system (\ref{4.9}), (\ref{4.90})
with condition (\ref{4.66}), a little problem is that these theorems are
valid only if Neumann boundary conditions equal zero at the entire lateral
boundary $S_{T}$ whereas in our case they are equal to zero only at $\Gamma
_{1,T}^{+}$, see (\ref{4.80}). On the other hand, Dirichlet boundary
conditions equal zero on the rest of the boundary $S_{T}\diagdown \Gamma
_{1,T}^{+}.$ Still, proofs of Carleman estimates of \cite{MFG4}, which are
actually taken from \cite{MFG1,MFG2}, can be straightforwardly rewritten for
the case when zero Neumann boundary conditions are given at $\Gamma
_{1,T}^{+}$ and zero Dirichlet boundary conditions are given at $%
S_{T}\diagdown \Gamma _{1,T}^{+}.$ Therefore, just as in the end of the
proof of Theorem 3.2, results of \cite{MFG4} again imply now that $%
\widetilde{u}\left( x,t\right) \equiv \widetilde{m}\left( x,t\right) \equiv
0 $ in $Q_{T}.$ \ \ $\square $

\section{Declarations}

\label{sec:6}

\textbf{Funding and/or Conflicts of interests/Competing interests.} The
author has no funding. The author does not have competing interests to
declare that are relevant to the content of this article.

\end{document}